 \documentclass[final,1p,times,review]{elsarticle}

\usepackage{lipsum}
\usepackage{amssymb}
\usepackage{graphicx}
\usepackage{epstopdf}
\usepackage{algorithmic}
\usepackage{caption,mathbbold,bbm,bm}
\usepackage{mathrsfs,amssymb,amsmath,amssymb,booktabs,array,xcolor,url,color,soul,multirow,multicol}
\usepackage{mathbbold,bbm}
\usepackage{slashbox}
\usepackage{stfloats}
\usepackage{amsfonts}
\usepackage{extarrows}
\usepackage{textcomp,balance}
\usepackage[justification=centering]{caption}
\usepackage{graphicx}
\usepackage{subfigure}
\usepackage[all]{xy}
\usepackage{chemarrow}

\newtheorem{theorem}{Theorem}[section]
\newtheorem{definition}[theorem]{Definition} 
 
\newtheorem{lemma}[theorem]{Lemma}
\newtheorem{corollary}[theorem]{Corollary}

\begin{document}
\begin{frontmatter}

\title{Persistence of Delayed Complex Balanced Chemical Reaction Networks\tnoteref{label1}}
\tnotetext[label1]{This work is supported by the National Nature Science Foundation of China under [Grant No. 11671418, 61611130124].}

 \author[label1]{Xiaoyu~Zhang}
 \ead{11735035@zju.edu.cn}
 \author[label2]{Chuanhou~Gao \corref{cor1}}
 \ead{gaochou@zju.edu.cn}
 \cortext[cor1]{Corresponding author}
\address[label1,label2]{School of Mathematical Sciences, Zhejiang University, Hangzhou 310027, China}



\begin{abstract}
In this paper, we derive two sufficient conditions to diagnose the persistence of two classes of delayed complex balanced chemical reaction network systems equipped with mass-action kinetics. One class is identified by $\dim Z_W=\dim\mathscr{S}-1$ ($Z_W$ is defined by \ref{eq:g} while $\mathscr{S}$ is the stoichiometric subspace of the network), the other class is identified through $\dim Z_W=0$ for any semilocking species set $W$ in the network. Then we also derive that delayed complex balanced systems with 2-d stoichiometric subspace are persistent. The results recur those proposed by Anderson et al. [D. F. Anderson, \textit{SIAM J. Appl. Math.}, 68 (2008), pp. 1464-1476; D. F. Anderson and A. Shiu, \textit{SIAM J. Appl. Math.}, 70 (2010), pp. 1840-1858] for checking the persistence of complex balanced reaction network systems without time delay. Further, we prove the above mentioned two classes of network systems are globally asymptotically stable at the corresponding positive equilibrium if the trajectory starts from a positive initial function. We illustrate the analysis by two numerical examples.    
\end{abstract}

\end{frontmatter}

\section{Introduction}
 Chemical reaction networks (CRNs) are popular mathematical models of several phenomena in population dynamics, system biology, epidemiology, telecommunications, and chemistry, etc. When equipped with mass-action kinetics \cite{Horn1972}, CRNs often give rise to a family of nonequilibrium polynomial ordinary differential equation (ODE) systems that capture the evolution of the concentrations of each species. It is usually impossible to find solutions of the dynamical equations. The past 40 years have witnessed the fast development of CRN theory since the pioneering work of Horn and Jackson \cite{Horn1972,Horn1974}, and Feinberg \cite{Feinberg1979,Feinberg1987}. One of the main goals is the characterization of the qualitative property of dynamical behaviors, for which to understand the long-time behavior of solutions is a central issue, such as the existence of stable equilibrium, persistence, periodic orbits, etc. In this paper, we are in line with the recent surge of interest in network persistence. 

Persistence is one of the important properties of the dynamics of CRNs. It means that if each species is present at the start of the reaction, no species will be completely used up in the course of the reaction. Mathematically, in the case of bounded trajectories, this property can be modeled as the requirement that there are no $\omega$-limit points in the boundary of the positive orthant for any trajectory starting in the interior of the positive orthant. The central importance of theory and practice has triggered an increasing interest in the study of persistence in recent decades. Feinberg \cite{Feinberg1987} even presented the well-known ``Persistence conjecture" that says \textit{any weakly reversible chemical reaction network is persistent}. A more special version is the ``Global attractor conjecture (GAC)", initiated by Feinberg and Horn \cite{Horn1974} and named GAC by Craciun et al. \cite{Craciun2009}, that states \textit{a complex balanced equilibrium contained in the interior of a positive stoichiometric compatibility class is a global attractor of the interior of that compatibility class}. If the former is true, so is the latter, since a complex balanced network must have a weakly reversible structure. There has been plenty of studies towards the proof of the above two conjectures, but they are still remaining open today. A foundation work came from Angeli et al. \cite{Angeli2007} that says if a boundary point of $\mathbb{R}^n_{\geq 0}$ (an element in $Q_W$, called a face of $\partial\mathbb{R}^{n}_{\geq 0}$) is an $\omega$-$limit$ point, then the species subset $W$ in which every species has zero concentration must be a semilocking set. This result is actually valid for networks with any structure. Anderson and his coworkers have made a systematic study on GAC \cite{Anderson2008,Anderson2010,Anderson2011}. They classified mass-action systems (MASs, CRNs equipped with mass-action kinetics) through the dimension of $\mathcal{P}_{x_0}\cap Q_{W}$, where $\mathcal{P}_{x_0}$ represents the stoichiometric compatibility class of the initial state $x_0$ and $Q_W$ means the face of $\partial\mathbb{R}^{n}_{\geq 0}$ in which the concentration of every species is zero if it belongs to the species subset $W$. When $\mathcal{P}_{x_0}\cap Q_{W}$ is empty, discrete, i.e., only containing some isolated vertexes ($\dim{(\mathcal{P}_{x_0}\cap Q_{W})}=0$) of $\mathcal{P}_{x_0}$ \cite{Anderson2008}, or a facet ($\dim{(\mathcal{P}_{x_0}\cap Q_{W})}=\dim\mathscr{S}-1$) of $\mathcal{P}_{x_0}$ \cite{Anderson2010} for any semilocking set $W$, a complex balanced MAS is persistent. Here, $\mathscr{S}$ is the stoichiometric subspace of the network. An immediate application of these results is that any complex balanced MAS with $\dim\mathscr{S}=2$ is persistent. This result is further extended to the case of $3$-dimensional complex balanced MASs by Pantea \cite{Pantea2012}. Also, Anderson \cite{Anderson2011} proved that any complex balanced MAS with a single linkage class is persistent. A thorough settlement of GAC was announced by Craciun \cite{Craciun2015} recently, but the proof is still under verification. 

The advance to prove these two conjectures are still underway. However, the studies on the persistence of CRNs are being made more extensively, enlarging to CRNs with non weakly reversible structure or time-delayed networks. Angeli et al. \cite{Angeli2007} reported that any MAS is persistent if there is a global conservative relation in the network and a local conservative relation among the species in $W$ for each semilocking set $W$. If the second condition is not true, they proposed a new concept ``dynamically non-emptiable", and fixed the result by adding one of the following conditions (a) all semilocking sets are dynamically non-emptiable; (b) there are no nested distinct locking sets. Using the same idea, Johnston and Siegel \cite{Johnson2011} extended this result by putting forward that MASs with all weak dynamical non-emptiable semilocking sets have persistence. Craciun et al. \cite{Craciun2013} defined a larger class of networks than weakly reversible ones, called endotactic structure, and Gopalkrishnan et al. \cite{G2014} further defined strongly endotactic network (a subclass of endotactic networks). They reported respectively that any endotactic MAS with two species and any strongly endotactic MAS are persistent. Based on their work, Zhang et al. \cite{Zhang19} proposed the notion of $\mathcal{W}_I$-endotactic network, which is a larger class than strongly endotactic network, and proved that any $1$-dimensional $\mathcal{W}_I$-endotactic network is persistent utilizing the Lyapunov function PDEs technique \cite{Fang2019}. 

Time delay in a reaction is another point taken into account in enriching the study of the persistence of CRNs. Certainly, the essential motivation is that time delay can simplify the dynamical equation of a very complicated network to a great extent and explain the dynamical process better. A typical example is a enzyme-substrate reaction network \cite{Hinch2004}, which takes place accompanying $N$ intermediate complexes $C_i$, i.e.,
 \begin{equation*}
S+E\xrightarrow{k_0}C_1\xrightarrow{k_1}\cdots\xrightarrow{k_{N-1}}C_N\xrightarrow{k_N}E+P,
\end{equation*}
where $S,E,P$ are substrate, enzyme and product, respectively. If the time delay is introduced, the original network can be modeled by a delayed enzyme-substrate reaction network in the form of
\begin{equation*}
S+E\xrightarrow{\rm{delay}}E+P
\end{equation*}
The network under study is simplified greatly. Astoundingly, time delay in the reaction will affect the dynamics in a very complicated means. It can cause a stable system unstable, even to generate oscillation, but it can be sometimes utilized to stabilize an unstable system \cite{Intro}. Therefore, time delay should be handled with great caution in modeling delayed mass-action systems (DeMASs). Liptk et al. \cite{G2018} studied the approximation of DeMASs based on the chain method \cite{Repin1965}. They further revealed the semistability of complex balanced MASs with arbitrary time delays \cite{G20182}. Komatsu and Nakajima \cite{H2019} made a pioneering work on the persistence of DeMASs. They derived a sufficient condition, i.e., the network is semiconservative with respect to every minimal semilocking set, to say a DeMAS persistent, which recurs the result presented by Angeli et al. \cite{Angeli2007, Angeli2011} for MASs without time delays.  

Despite the above facts, there still need more concerns on DeMASs, especially in persistence analysis. In this paper, we manage to study the persistence of two classes of delayed complex balanced MASs (DeCBMASs), identified by $\dim Z_W=\dim\mathscr{S}-1$ ($Z_W$ is defined by \ref{eq:g}) and $\dim Z_W=0$, respectively. We firstly present the $\omega$-$limit$ set theorem for DeMASs. Then two sufficient conditions are given to diagnose the persistence of the mentioned two classes of DeCBMASs. These results are natural analogous of those suggested by Anderson et al. \cite{Anderson2008, Anderson2010} for MASs without time delays. The reminder of this paper is organized as follows. \ref{sec2} recalls some basic concepts and results of CRNs, MASs, and DeMASs. This is followed by the $\omega$-$limit$ set theorem for DeMASs in \ref{sec3}. In \ref{sec4}, two classes of DeCBMASs, labeled by different dimensions of $Z_W$, are proved to have persistence. \ref{sec5} illustrates the results through two numerical examples. Finally, we conclude the paper in \ref{sec6}.

~\\~\\
\noindent{\textbf{Mathematical Notation:}}\\
\rule[1ex]{\columnwidth}{0.8pt}
\begin{description}
\item[\hspace{+0.8em}{$\mathbb{R}^n, \mathbb{R}^n_{\geq 0},\mathbb{R}^n_{>0}$}]: $n$-dimensional real space, non-negative and positive real space, respectively.
\item[\hspace{+0.8em}{$\mathbb{Z}^n_{\geq 0}$}]: $n$-dimensional non-negative integer space.
\item[\hspace{+0.8em}{$\bar{\mathscr{C}}_{+}=C([-\tau,0];\mathbb{R}^{n}_{\geq 0}), \mathscr{C}_{+}=C([-\tau,0];\mathbb{R}^{n}_{>0})$}]: the non-negative, positive continuous function vectors defined on the interval $[-\tau,0]$, respectively.
\item[\hspace{+0.8em}{$x^{y_{\cdot i}}$}]: $x^{y_{\cdot i}}\triangleq\prod_{j=1}^{n}x_{j}^{y_{ji}}$, where $x,y_{\cdot i}\in\mathbb{R}^{n}$.
	\item[\hspace{+0.8em}{$\mathrm{Ln}(x)$}]: $\mathrm{Ln}(x)\triangleq\left(\ln{x_{1}}, \cdots, \ln{{x}_{n}} \right)^{\top}$, where $x\in\mathbb{R}^{n}_{>{0}}$.
	\item[\hspace{+0.8em}{supp$~x$}]: Support set, defined by $\text{supp}~ x=\{i|x_{i}\neq 0\},~\forall x\in \mathbb{R}^{n}_{\geq 0}$.
\item[\hspace{+0.8em}{$\mathbbold{0}_{n}$}]: An $n$-dimensional vector with each element to be zero.
\item[\hspace{+0.8em}{$0^{0}$}]: The result is defined by $1$.
\end{description}
\rule[1ex]{\columnwidth}{0.8pt}
\section{Preliminaries}\label{sec2}
In this section, we formally introduce CRN and the associated MASs and DeMASs \cite{Feinberg1979,Anderson2010}.
\subsection{CRNs and MASs}
A CRN is a list of chemical reactions taking place among constituent chemical species, and is defined as follows.
\begin{definition}[CRN] A CRN is given by a triple of finite sets $(\mathcal{S,C,R})$ where 
\begin{itemize}
    \item $\mathcal{S}=\bigcup_{j=1}^{n}\{X_{j}\}$ is a set of species that participate reactions;
    \item $\mathcal{C}=\bigcup_{i=1}^r\{y_{\cdot i},y'_{\cdot i}\}$ is a set of linear combinations of species with non-negative integer coefficients, referred to as complexes, with $y_{\cdot i},y'_{\cdot i}\in \mathbb{Z}^n_{\geq 0}$ appearing on the left and right side of the reaction arrow, respectively, and satisfying $$\bigcup_{i=1,\cdots,r}\{\mathrm{supp} ~{y_{\cdot i}}\cup \mathrm{supp}~{y'_{\cdot i}}\}=\mathcal{S};$$
\item $\mathcal{R}=\bigcup_{i=1}^{r}\{y_{\cdot i}\rightarrow y'_{\cdot i}\}$ is a set of reactions with $\forall~i,~y_{\cdot i}\neq y'_{\cdot i}$ and $\forall~y_{\cdot i}\in\mathcal{C},~\exists~ y_{\cdot j}\in\mathcal{C}$ supporting either $y_{\cdot i}\rightarrow y_{\cdot j}\in\mathcal{R}$ or $y_{\cdot j}\rightarrow y_{\cdot i}\in\mathcal{R}$. 
\end{itemize}
\end{definition}
The definition suggests the $i$th reaction in $(\mathcal{S,C,R})$ to be written as:
$$y_{1i}X_{1}+\cdots+y_{ni}X_{n}\rightarrow y'_{1i}X_{1}+\cdots+y'_{ni}X_{n},$$
where $y_{ji},~y'_{ji}$ are the coefficients of $X_j$ in the source complex $y_{\cdot i}$ and product complex $y'_{\cdot i}$, respectively.

A CRN $(\mathcal{S,C,R})$ can be associated with a directed graph by considering the set of complexes as nodes and the set of reactions as directed edges, often referred to as a reaction graph. They are a one-to-one correspondence. Each connected component in the reaction graph is termed a linkage class of the network. Two kinds of special networks may be defined based on the structure of the reaction graphs, which we are concerned about in this paper.

\begin{definition}[weakly reversible and reversible CRN]\label{def:2reversible} A CRN $(\mathcal{S,C,R})$ is called weakly reversible if its reaction graph is strongly connected. Namely, for any reaction $y_{\cdot i}\to y'_{\cdot i}\in\mathcal{R}$ there exists a chain of reactions starting from $y'_{\cdot i}$ and ending with $y_{\cdot i}$, i.e., $y'_{\cdot i} \to y_{\cdot i_{1}}\in\mathcal{R}$, $\cdots$, $y_{\cdot i_{m}}\to y_{\cdot i}\in\mathcal{R},~m\leq r$. In particular, if $\forall~y_{\cdot i}\to y'_{\cdot i}\in\mathcal{R}$, $\exists~y'_{\cdot i}\to y_{\cdot i}\in\mathcal{R}$, then this network is called reversible.
\end{definition}

Clearly, reversible CRNs are a subclass of weakly reversible ones. We further define the reaction vector by $y'_{\cdot i}-y_{\cdot i}$ and the stoichiometric subspace by \begin{equation*}
    \mathscr{S}=\mathrm{span}\{y'_{\cdot i}-y_{\cdot i}:~i=1,\cdots,r\}.
\end{equation*}

The popular mathematical modeling for CRNs is to capture the evolution of the concentrations of each species, labeled by $x\in\mathbb{R}^n_{\geq 0}$. By defining a real-valued function $R:\mathbb{R}^n_{\geq 0}\to \mathbb{R}_{\geq 0}$ to evaluate the reaction rates and further based on the balance laws, we can write the dynamical equation of a CRN to be
\begin{equation}\label{eq:de}
    \dot{x}(t)=\sum_{i=1}^{r}R_i(x(t))(y'_{\cdot i}-y_{\cdot i}).
\end{equation}
The model is essentially a coupled set of some nonlinear ordinary differential equations (ODEs), and its integral form takes
\begin{equation}\label{eq:ifde}
    x(t)=x_0+\sum_{i=1}^{r}\left(\int^{t}_{0}R_i(x(s))ds\right)(y'_{\cdot i}-y_{\cdot i}),
\end{equation}
where $x_0=x(0)\in\mathbb{R}^n_{\geq 0}$ represents the initial state. The integral form of \ref{eq:ifde} implies $x(t)-x_0\in \mathscr{S}$, i.e., the trajectories of the network system do not run all through the whole $n-$dimensional nonnegative real space $\mathbb{R}^n_{\geq 0}$, but only within an invariant set.  

\begin{definition}[stoichiometric compatibility class]\label{def:SCC} For a CRN $(\mathcal{S,C,R})$, let $x_0\in\mathbb{R}^{n}_{\geq 0}$, the sets $$\mathcal{P}_{x_0}:=\{x\in\mathbb{R}^{n}\vert x-x_{0}\in \mathscr{S}\},$$ $\bar{\mathcal{P}}_{x_0}^+:=\mathcal{P}_{x_0}\cap \mathbb{R}^{n}_{\geq 0}$ and $\mathcal{P}_{x_0}^+:=\mathcal{P}_{x_0}\cap \mathbb{R}^{n}_{>0}$ are named the stoichiometric compatibility class, the nonnegative stoichiometric compatibility class, and the positive stoichiometric compatibility class of $x_0$, respectively. 
\end{definition}

It is clear that the trajectory of the network system \ref{eq:ifde} will evolve in the invariant set $\bar{\mathcal{P}}_{x_0}^+$ if it starts from $x_0$. The stoichiometric compatibility class $\mathcal{P}_{x_0}$ is indeed a polyhedron. As done in \cite{Anderson2010}, we define the dimension of a polyhedron in $\mathbb{R}^{n}$ by the dimension of the span of the translate of this polyhedron which contains the origin, i.e., $\dim\mathcal{P}_{x_{0}}=\dim\mathscr{S}$. 

The widely-used kinetics for measuring the reaction rates is the mass-action kinetics, which assigns each reaction $y_{\cdot i}\to y'_{\cdot i}$ a positive real number $k_i$, called a reaction rate constant, and quantifies the reaction rate according to the power law,
\begin{equation}\label{eq:mas}
    R_i(x)=k_ix_1^{y_{1i}}x_2^{y_{2i}}\cdots x_n^{y_{ni}}=:k_ix^{y_{.i}}.
\end{equation}
We call CRNs equipped with mass-action kinetics to be MASs, referred to as 
$(\mathcal{S,C,R},$
$k)$ in the context. By inserting \ref{eq:mas} into \ref{eq:de}, we get the dynamics of a MAS
\begin{equation}\label{eq:mde}
    \dot{x}(t)=\sum^{r}_{i=1}k_ix(t)^{y_{.i}}(y'_{\cdot i}-y_{\cdot i}). 
\end{equation}

\begin{definition}[equilibrium and complex balanced equilibrium]\label{def:equilibrium}
A concentration vector $\bar{x}\in\mathbb{R}^n_{\geq 0}$ is called an equilibrium for the dynamics $\dot{x}(t)=\sum^{r}_{i=1}k_ix(t)^{y_{.i}}(y'_{\cdot i}$
$-y_{\cdot i})$ if $$\sum^{r}_{i=1}k_i\bar{x}^{y_{.i}}(y'_{\cdot i}-y_{\cdot i})=0,$$ and is a complex balanced equilibrium if $\forall~z\in\mathcal{C}$, there exists
\begin{equation}\label{eq:cb}
    \sum_{i:~y_{\cdot i}=z}k_i(\bar{x})^z=\sum_{i:~y'_{\cdot i}=z}k_i(\bar{x})^{y_{.i}}.
\end{equation}
\end{definition}
Feinberg \cite{Feinberg1979} revealed that if one equilibrium in a MAS is complex balanced, so are others. 
A MAS having a complex balanced equilibrium is called a complex balanced MAS. For the sake of distinction, we call $\bar{x}\in\mathbb{R}^n_{> 0}$ a positive (complex balanced) equilibrium, and others the boundary (complex balanced) equilibria. It has been proven that each positive complex balanced equilibrium is locally asymptotically stable relative to its stoichiometric compatibility class using the well-known pseudo-Helmholtz function 
\begin{equation}
    H(x)=\sum_{j=1}^{n}\left(\bar{x}_{j}-x_{j}-x_j\ln{\frac{\bar{x}_{j}}{x_j}}\right)
\end{equation}
as a Lyapunov function \cite{Horn1972}. 

We define the facet of a polyhedron to be better associated with the boundary equilibria in $(\mathcal{S,C,R},k)$. 

\begin{definition}[facet]
Let $Q$ be a polyhedron, a facet of $Q$ is the face of this polyhedron with dimension one less than it. Further denote $W\subseteq \mathcal{S}$ and $W\neq \varnothing$, where $\mathcal{S}$ is the species set in $(\mathcal{S,C,R},k)$, the set $Q_W$ is defined by
\begin{equation}\label{def:Qw}
    Q_W=\{x\in \mathbb{R}^{n}\vert x_i=0,~ X_i\in W;~x_i\neq 0,~X_i\notin W\}.
\end{equation}  
\end{definition}
Utilizing $Q_W$, the boundary of $\mathbb{R}^{n}_{\geq 0}$ may be expressed as
$$\partial\mathbb{R}^{n}=\cup_{W\subseteq\mathcal{S}}Q_{W}.$$ 
and the face of the stoichiometric compatibility class of $x_0$ as $Q_W\cap \mathcal{P}_{x_0}$. Clearly, each face of $\mathcal{P}_{x_0}$ is also a polyhedron with a lower dimension. 

The following concept plays an important role in the subsequent persistence analysis. 

\begin{definition}[semilocking set and locking set]
For a CRN $(\mathcal{S,C,R})$, a non-\\empty symbol set $W\subseteq\mathcal{S}$ is called a semilocking set if it satisfies: $W\cap\mathrm{supp}~y_{\cdot i}\neq \varnothing$ if $W\cap\mathrm{supp}~y'_{\cdot i}\neq\varnothing$.
 If for any reaction $y_{\cdot i}\to y'_{\cdot i}$, there is $W\cap\mathrm{supp}~y_{\cdot i}\neq \varnothing$, $W$ is called a locking set.  
\end{definition}

 \subsection{DeMASs}
There usually exists a time delay in chemical reactions, i.e., the production occurs later than the consumption. This factor will affect the dynamics of a MAS but have no effect on the network structure. Therefore, a DeMAS has some same concepts, like the reversibility, weak reversibility, stoichiometric subspace, linkage class, etc. as the MAS without time delay. Assume $\tau_i\geq 0$ to represent the time delay of the $i$th reaction in a MAS $(\mathcal{S,C,R},k)$, then its dynamical equation is in the form of
 \begin{equation}\label{eq:dde}
     \dot{x}(t)=\sum^{r}_{i=1}k_i
     [(x(t-\tau_i))^{y_{\cdot i}}y'_{\cdot i}-(x(t))^{y_{\cdot i}}y_{\cdot i}],~~~t\geq 0.
 \end{equation}
 For convenience, we use a five-tuple $(\mathcal{S,C,R},k,\tau)$ to denote a DeMAS where $\tau\in\mathbb{R}^{r}$. The dynamics of \ref{eq:dde} will reduce to \ref{eq:mde} if all the reactions delay are equal to zero. Obviously, the solution to the ODE \ref{eq:dde} relies on the initial function $x(s)=\psi(s)$ where $s\in[-\tau_{\text{max}},0]$ and $\tau_{\text{max}}=\max_{1\leq i\leq r}\tau_{i}$. Hence, the solution space of \ref{eq:dde} is a functional space \cite{G2018,G20182}, given by $\bar{\mathscr{C}}_{+}:=C([-\tau_{\text{max}},0];\mathbb{R}^{n}_{\geq 0})$. Also, denote $\mathscr{C}_{+}:=C([-\tau_{\text{max}},0];\mathbb{R}^{n}_{> 0})$ and $\mathscr{C}:=C([-\tau_{\text{max}},0];\mathbb{R}^{n})$, and define the norm in $\mathscr{C}$ by
  \begin{equation}
     \Vert\psi\Vert:=\sup_{-\tau_{\text{max}}\leq s \leq 0}\vert\psi(s)\vert
 \end{equation}
 with $\psi\in\mathscr{C}$. We further define some concepts related to dynamics for DeMASs that are different from those for MASs. 
  \begin{definition}[delayed stoichiometric compatibility class] For a DeMAS $(\mathcal{S,C,}$
  $\mathcal{R},k,\tau)$, let $\psi\in\bar{\mathscr{C}_{+}}$. Then the stoichiometric compatibility class (respectively, positive stoichiometric compatibility class) of $\psi$ is defined by
 \begin{equation}\label{eq:Dpsi}
 \mathcal{D}_{\psi}:=\{\phi\in \bar{\mathscr{C}}_{+} \vert c_{a}(\phi)=c_{a}(\psi),~\forall ~a\in \mathscr{S}^{\bot}\},
 \end{equation} 
 $$(respectively,~\mathcal{D}^+_{\psi}:=\{\phi\in \mathscr{C}_{+}\vert c_{a}(\phi)=c_{a}(\psi),~\forall ~a\in \mathscr{S}^{\bot}\} )$$
 where $\mathscr{S}^{\bot}$ is the orthogonal complement of $\mathscr{S}$, and the functional $c_{a}:\bar{\mathscr{C}_{+}}\rightarrow \mathbb{R}$ (respectively, $c_{a}:\mathscr{C}_{+}\rightarrow \mathbb{R}$) is defined by
 \begin{equation*}
     c_{a}(\psi):=a^{\top}\left[\psi(0)+\sum_{y_{\cdot i}\rightarrow y'_{\cdot i}}\left(\int_{-\tau_{\text{max}}}^{0}k_i\psi(s)^{y_{\cdot i}}ds\right)y_{\cdot i}\right].
 \end{equation*}
 \end{definition}
 
 The equilibrium and complex balanced equilibrium in a DeMAS share the same definitions as those in a MAS \cite{H2019}, i.e., following \ref{def:equilibrium}. We call a DeMAS having a complex balanced equilibrium a DeCBMAS. For a DeCBMAS, each stoichiometric compatibility class contains a unique positive equilibrium \cite{G20182}. Moreover, Liptk et al. \cite{G20182} proved that each positive complex balanced equilibrium in a DeCBMAS is also locally asymptotically stable. The used Lyapunov function $V:\mathscr{C}_{+}\rightarrow \mathbb{R}_{\geq 0}$ is given by  
 \begin{equation}\label{eq:Vd}
\begin{split}
V(\psi)&=\sum_{j=1}^{n}(\psi_{j}(0)(\ln(\psi_{j}(0))-\ln(\bar{x}_{j})-1)+\bar{x}_{j})\\
&+\sum_{i=1}^{r}k_{i}\int^{0}_{-\tau_{i}}\left\{(\psi (s))^{y_{\cdot i}}[\ln((\psi(s)^{y_{\cdot i}})-\ln(\bar{x}^{y_{\cdot i}})-1]+\bar{x}^{y_{\cdot i}}\right\} ds,
\end{split}
\end{equation}
named the delayed pseudo-Helmholtz function. 

Finally, we define persistence for DeMASs. 

\begin{definition}[persistence]
A DeMAS $(\mathcal{S,C,R},k,\tau)$ described by \ref{eq:dde} has persistence if any forward trajectory $x^{\psi}(t)\in \mathbb{R}^{n}_{\geq 0}$ with a positive initial condition $\psi\in \mathscr{C}_{+}$ satisfies
\begin{equation*}
\liminf_{t\rightarrow \infty}x^{\psi}_{j}(t)>0~~~for~all~j\in\{1,\cdots,n\}.
\end{equation*}
\end{definition}
In the case of bounded systems, the persistence can be characterized through the $\omega$-$limit$ set. 
\begin{definition}[$\omega$-$limit$ set]
The $\omega$-$limit$ set for the trajectory $x^{\psi}(t)$ with a positive initial condition $\psi\in \mathscr{C}_{+}$ is
\begin{equation*}
    \omega(\psi):=\{\phi\in \bar{\mathscr{C}}_{+}~|~x^{\psi}_{t_{N}}\rightarrow \phi,~\text{for some sequence} ~t_{N}\rightarrow \infty~with ~t_{N}\in\mathbb{R}\}.
\end{equation*}
Each element in the $\omega$-$limit$ set is an $\omega$-$limit$ point of this trajectory.
\end{definition}

\begin{definition}[persistence for bounded trajectories]
For a DeMAS $(\mathcal{S,C,R},$
$k,
\tau)$ with bounded trajectories, it is persistent if 
\begin{equation}
 \omega(\psi)\cap(\cup_{W} L_W)=\varnothing
 , ~~~~\forall~\psi\in\mathscr{C}_{+},
\end{equation} 
where 
\begin{equation}\label{def:Lw}
   L_W=\left\{w\in \bar{\mathscr{C}}_{+}\big|\substack{w_j(s)=0,~ X_j\in W,\\w_j(s)\neq 0,~X_j\notin W,}~~\forall s\in [-\tau_{\rm{max}},0]\right\}.
\end{equation}
\end{definition}

\section{$\omega$-$limit$ set theorem of DeMASs}\label{sec3} 
The $\omega$-$limit$ set is an important tool in diagnosing persistence of a MAS with boundary trajectories. Siegel and Maclean \cite{Siegel2000} presented the $\omega$-$limit$ set theorem that says the $\omega$-$limit$ set, denoted by $\omega(x_0)$, for a complex balanced MAS, is either a set of the boundary equilibria or the single positive equilibrium. Anderson et al. \cite{Anderson2008,Angeli2007} further proved that the $\omega$-limit point, i.e., any element in $\omega(x_0)$, can only exist in the face of the stoichiometric compatibility class of $x_0$, i.e., $Q_W\cap \mathcal{P}_{x_0}$, when $W$ is a semilocking set. This result also applies to DeMASs \cite{H2019}, i.e., the $\omega$-limit point in $\omega(\psi)$ can only exist in $L_W\cap\mathcal{D}_{\psi}$ in the case that $W$ is a semilocking set. It is expected that the $\omega$-$limit$ set theorem also holds for DeCBMASs.

We begin with a general result (not limited to MASs) for obtaining the $\omega$-$limit$ set theorem for DeCBMASs.   
\begin{lemma}\label{lm:con}
Consider a delayed dynamical system given by \ref{eq:dde}. Let  $\varGamma\in \mathcal{F}$ represent any trajectory of the system starting from an initial function $\psi\in \mathcal{F}$, where $\mathcal{F}$ is an open subset of $\bar{\mathscr{C}}_{+}$. Then the $\omega$-limit set $\omega(\psi)$ for the trajectory $\varGamma$ is a closed set. Further, if $\varGamma$ is contained in a compact subset of $\mathcal{F}$, then $\omega(\psi)$ is a non-empty, connected and compact subset of $\mathcal{F}$.
\end{lemma}
\noindent{\textbf{Proof.}}
 If a function sequence $\omega_m\in \omega(\psi)$ has $\omega_m\rightarrow \omega$, for each $m$, we can find a time sequence $t^{(m)}_{i}$ such that:
\begin{equation*}
    \lim_{i\rightarrow \infty}\Gamma(t_{i}^{(m)},\psi)=\omega_m.
\end{equation*}
In other words, for each $m$, there exists $K(m)$ such that when $i\geq K(m)$, 
\begin{equation*}
\sup_{-\tau_{max}\leq s \leq 0}\vert \Gamma(t_{i}^{(m)}-s,\psi)-\omega_{m}(s)\vert<1/m.
\end{equation*}
Let $t_m=t_{K(m)}^{(m)}$, we have:
\begin{align*}
     & \sup_{-\tau_{max}\leq s \leq 0}\vert \Gamma(t_m-s,\psi)-\omega(s)\vert \leq\\ &\sup_{-\tau_{max}\leq s\leq 0}(\vert\Gamma(t_m-s,\psi)-\omega_m(s)\vert+\vert\omega_m(s)-\omega(s)\vert)
      \leq 1/m+\Vert \omega_m-\omega\Vert
\end{align*}
  Above equation trends to zero when $m$ trends to infinity. We can easily conclude that $\omega\in \omega(\psi)$. Therefore, $\omega(\psi)$ is a closed set.
  If $\Gamma$ is contained in a compact subset $K$ and we can choose a sequence $t_m$ such that $\Gamma(t_m,\psi)\rightarrow \omega\in \omega(\psi)$. $\omega(\psi)\subset K$ and is compact. 
  
 Now we assume that $\omega(\psi)$ is not connected. In this case two nonempty, closed sets $A$,$B$ can be found such that $\omega(\psi)=A\cup B$ and $\delta$ is a positive finite distance between $A$ and $B$ where 
  \begin{equation*}\label{d}
      \delta(A,B)=\inf_{\varphi_1\in A,\varphi_2\in B}\sup_{s\in [-\tau_{max},0]}\vert \varphi_1(s)-\varphi_2(s)\vert
 \end{equation*}
  Since $A$ and $B$ are both the subsets of $\omega(\psi)$, we can find a big enough $t_{m_{1}}$ such that $\Gamma(t_{m_{1}})$  within the $\delta(A,B)/2$-neighbourhood of A, and a big enough $t_{m_2}$ such that $\Gamma (t_{m_{2}})$ out of $\delta(A,B)/2$-neighbourhood of A. From the continuity of $\delta(\Gamma(t),A)$, we can get there exists a $t_{m}$ lead to $\delta(\Gamma(t_m),A)=\delta/2$. Thus let $\lim_{t_m\rightarrow+\infty}\Gamma(t_m)=\omega$, $\omega\notin A\cup B$ is also an $\omega$-$limit$ point of $\Gamma(t,\psi)$. This contradicts to $\omega(\psi)=A\cup B$. So the connectivity of $\omega(\psi)$ is proved.$\Box$
  
Based on \ref{lm:con}, we then could reach the $\omega$-$limit$ set theorem for DeCBMASs.

\begin{theorem}[$\omega$-limit set theorem for DeCBMASs]\label{thm:wl}
The $\omega$-limit set $\omega(\psi)$ of a DeCBMAS is either a set of the constant boundary equilibria or the unique constant positive complex balanced equilibrium relative to the stoichiometric compatibility class $\mathcal{D}_{\psi}$, where $\psi$ is the initial data of the trajectory. 
\end{theorem}  

\noindent{\textbf{Proof.}}
From the result in \cite{G2018}, we know that $V:\mathscr{C}_{+}\rightarrow \mathbb{R}_{+}$ in \ref{eq:Vd} is a convex function and decrease along the trajectory. So the trajectory must be bounded as the same reason with the non-delayed system.
\cite{G2018} tells us that
only equilibrium can be the $\omega$-$limit$ point in positive stoichiometric compatibility class. Now we consider the points on the boundary in two cases:

Firstly, if $W$ is a locking set, according to the definition of the boundary $L_W$, all the functions in $L_W$ are all equilibria.
And it is obvious that $\dot{V}(\phi)=0$ for any $\phi\in L_W$. Now we prove that $\omega$-$limit$ set can only be constant equilibria. If $\phi\in \omega(\psi)\cap L_W$, we can find a $\phi_1\in \omega(\psi)\cap L_W$ such that $x^{\phi}(t)=x^{\phi_1}(t-\tau_{max})$.  Now we consider the dynamics of the $x^{\phi_1}(t)$ for any $t\in [0,\tau_{max}]$
\begin{equation}\label{eq:s}
   \dot{x}^{\phi_1}(t)=\sum^{r}_{i=1}k_i\left(x(t-\tau_{y_{\cdot i}\rightarrow y'_{\cdot i}})\right)^{y_{\cdot i}}y'_{\cdot i}-\sum^{r}_{i=1}k_{i}\left(x(t)\right)^{y_{\cdot i}}y_{\cdot i}=0,~~\forall~t\in~[0,\tau_{max}].
\end{equation}
Above equation holds as $x(t-\tau)$ and $x(t)$ are all in $L_W$ and all the $y_{\cdot i}$ satisfy $\mathrm{supp}~y_{\cdot i}\cap W\neq\emptyset$ and \ref{eq:s} reveals that $\dot{\phi}\equiv 0$. Therefore, if $\phi\in \omega(\psi)\cap L_W$, $\phi$ must be a constant equilibrium.

The second part contributes to prove the case that $W$ is a semilocking set but not a locking set. In this case, we can find a subnetwork $\mathcal{N}_W=\mathcal(S',C',R')$ where
$$\mathcal{R'}=\{y\rightarrow y'\vert \text{supp}~y\cap W= \emptyset,~\text{supp}~y'\cap W= \emptyset\}$$ with  $\vert\mathcal{R'}\vert=r'$. $\mathcal{S'}$, $\mathcal{C'}$ are the set of species and complex appear in the reactions in $\mathcal{R'}$ respectively. Without loss of generality, we can assume $\mathcal{S'}=\{X_1,...,X_{N'}\}$. This $\mathcal{N}_W$ is obviously a complex balanced network and its dynamical function is \ref{eq:de}.
From the properties of the complex balanced delayed system, each stoichiometric compatibility class only has one positive constant equilibrium.
From the Lyapunov function of $\mathcal{N}_W$ --- $V_1$, defined as:
\begin{align*}
V_1(\phi')&=\sum_{j=1}^{N'}(\phi'_{j}(0)(\ln(\phi'_{j}(0))-\ln(\bar{x}_{j})-1)+\bar{x}_{j})\\
&+\sum_{i=1}^{r'}k_{i}\int^{0}_{-\tau_{i}}\{(\phi' (s))^{y_{\cdot i}}[\ln((\phi'(s))^{y_{\cdot i}})-\ln(\bar{x}^{y_{\cdot i}})-1]+\bar{x}^{y_{\cdot i}}\} ds\leq 0
\end{align*}
form \cite{G2018} we know above equation with equality if and only if for each $i=1,...,r'$, 
\begin{equation*}
  \left(\frac{\phi'(0)}{\bar{\phi}'}\right)^{y'_{\cdot i}}=\left(\frac{\phi'(-\tau_i)}{\bar{\phi}'}\right)^{y_{\cdot i}}
 \end{equation*}
Then by using above equation and rewrite \ref{eq:de} into following form we can obtain that:
\begin{align*}
      \dot{x}(t)&=\sum_{c\in \mathcal{C'}}\left[ \sum_{c=y'_{\cdot i}}k_i\bar{x}^{y_{\cdot i}}\left(\frac{x(t-\tau_i)}{\bar{x}}\right)^{y_{\cdot i}}-\sum_{c=y_{\cdot i}}k_i\bar{x}^{y_{\cdot i}}\left(\frac{x(t)}{\bar{x}}\right)^{y_{\cdot i}}\right]c\\
      &=\sum_{c\in \mathcal{C'}}\left(\frac{x(t)}{\bar{x}}\right)^{c}\left[\sum_{c=y'_{\cdot i}}k_i\bar{x}^{y_{\cdot i}}-\sum_{c=y_{\cdot i}}k_i\bar{x}^{y_{\cdot i}}\right]c=0
\end{align*}
  The last equation equal to zero can be derived from the  
definition of the complex balanced equilibrium. 
From above discussion, we can obtain that only constant equilibrium can be the $\omega$-$limit$ function in $L_W$ as $\dot{V}_1(\phi)$ is zero if and only if $\phi$ is the unique positive equilibrium of $\mathcal{N}_W$. We can easily see $\dot{V}_1(\phi)=\lim\limits_{\gamma \rightarrow \varphi}\dot{V}(\gamma)$  where $\gamma\in \mathscr{C}_+$, $\varphi\in \bar{\mathscr{C}_+}\cap L_W$ and $\phi$ is the function $\varphi $ restricted on $W^{c}$. Hence, from the continuity of $V$, we obtain that $\dot{V}(\varphi)=\lim\limits_{\gamma\rightarrow \varphi}\dot{V}(\gamma)=\dot{V}_1(\phi)=0$ if and only if $\varphi$ is the constant equilibrium in $L_W$. And the non-equilibrium function $\varphi_1$ in $L_W$ can not be the $\omega$-$limit$ function of any trajectories due to $\dot{V}(\varphi_1)<0$. Thus we complete our proof.
$\Box$

The above result is analogous to that given in \cite{Siegel2000} for complex balanced MASs. 
\section{Persistence of two classes of DeCBMASs}\label{sec4}
In this section, we will derive two sufficient conditions for diagnosing the persistence of DeCBMASs, which might be thought of as the parallel results given by Anderson et al. \cite{Anderson2008,Anderson2010} for complex balanced MASs. They said that for a complex balanced MAS $(\mathcal{S,C,R},k)$ with $x_0$ as the initial point, it is persistent if one of the following conditions holds in the case of any nonempty $W\subseteq\mathcal{S}$ being a semilocking set: 1) $Q_W\cap\mathcal{P}_{x_0}=\varnothing$; 2) $Q_W\cap\mathcal{P}_{x_0}$ is a facet of $\mathcal{P}_{x_0}$; 3) $Q_W\cap\mathcal{P}_{x_0}$ is discrete (namely, being vertexes).  

We define the following set to sort DeCBMASs.
\begin{definition}
Consider a DeCBMAS $(\mathcal{S,C,R},k,\tau)$ described by \ref{eq:dde} with an initial function $\psi\in\bar{\mathscr{C}_{+}}$. Define a vector function by
$$g(x)=x(0)+\sum_{i=1}^{r}k_i\left(\int_{-\tau_i}^{0}(x(s))^{y_i}ds\right)y_{\cdot i},$$ and further define $Z_W$ by
\begin{equation}\label{eq:g}
    Z_W=\left\{v\vert v=g(x_1(s))-g(x_2(s)),x_1(s), x_2(s)\in L_W\cap\mathcal{D}_{\psi}, s\in[-\tau_{\rm{max}},0]\right\}.
\end{equation}
If $\dim{Z_W}=\dim\mathscr{S}-1$, then $L_W\cap\mathcal{D}_{\psi}$ is called a ``facet" of the stoichiometric compatibility class $\mathcal{D}_{\psi}$. If $\dim{Z_W}=0$, then $L_W\cap\mathcal{D}_{\psi}$ is called a ``vertex" of the stoichiometric compatibility class $\mathcal{D}_{\psi}$.
\end{definition}

In the following, we will use $Z_W$ to label some DeCBMASs that have persistence.

\subsection{DeCBMASs with $\dim{Z_W}=\dim{\mathscr{S}-1}$ for any semilocking set $W\subseteq \mathcal{S}$} 
We firstly consider the DeCBMASs in which $L_W\cap\mathcal{D}_{\psi}$ is a facet of the stoichiometric compatibility class $\mathcal{D}_{\psi}$. 

\begin{lemma}
For a DeMAS $(\mathcal{S,C,R},k,\tau)$ with the dynamics of \ref{eq:dde}, let $W\subseteq \mathcal{S}$ and $W\neq\varnothing$. If $L_W\cap\mathcal{D}_{x_0(s)}$ is a facet of stoichiometric compatibility class $\mathcal{D}_{x_0(s)}$ in $(\mathcal{S,C,R},k,\tau)$,  $Q_W\cap \mathcal{P}_{x_0}$ given in \ref{def:Qw} is also a facet of the stoichiometric compatibility class $\mathcal{P}_{x_0}$ of $x_0$ in the corresponding MAS $(\mathcal{S,C,R},k)$ without time delay.
\end{lemma}
\noindent{\textbf{Proof.}}
 We assume that the set $\{w_1,w_2,...,w_{s-1}\}$ where any $w_i$ and $w_j$ are independent can denote the elements of $Z_W$ by linear combination. And without loss of generality, we assume $W=\{X_{d+1},...,X_{n}\}$. For any $x_0(s)$ and $x_1(s)$ in $L_W\cap\mathcal{D}_{x_0(s)}$, we have:
\begin{align*}
g&(x_0)-g(x_1)=x_0(0)-x_1(0)+\sum^{r}_{i=1}k_i\int^{0}_{-\tau_i}(x_0(s))^{y_{\cdot i}}ds y_{\cdot i}-\sum^{r}_{i=1}k_i\int^{0}_{-\tau_i}(x_1(s))^{y_i}dsy_{\cdot i}\\
=&\dbinom{\tilde{x}_0(0)}{\mathbbold{0}_{n-d}}-\dbinom{\tilde{x}_1(0)}{\mathbbold{0}_{n-d}}+\sum_{\mathrm{supp}~y_{\cdot i}\cap W=\emptyset}k_{\cdot i}\int^{0}_{-\tau_i}\left(\left(\tilde{x}_0\left(s\right)\right)^{\tilde{y}_{\cdot i}}-\left(\tilde{x}_1\left(s\right)\right)^{\tilde{y}_{\cdot i}}\right)ds\dbinom{\tilde{y}_{\cdot i}}{\mathbbold{0}_{n-d}}
\end{align*}
Let $(\tilde{x}_0(s))^{\tilde{y}_{\cdot i}}$ smaller than $(\tilde{x}_1(s))^{\tilde{y}_{\cdot i}}$, otherwise we can exchange $x_0$ and $x_1$. Let
$$x_2=\dbinom{\tilde{x}_1(0)}{\mathbbold{0}_{n-d}}-\sum_{\mathrm{supp}~y_{\cdot i}\cap W=\emptyset}k_i\int^{0}_{-\tau_i}\left(\left(\tilde{x}_0\left(s\right)\right)^{\tilde{y}_{\cdot i}}-\left(\tilde{x}_1\left(s\right)\right)^{\tilde{y}_{\cdot i}}\right)ds\dbinom{\tilde{y}_{\cdot i}}{\mathbbold{0}_{n-d}}.$$ So we can easily verify that $x_2$ is a point in $Q_W\cap \mathcal{P}_{x_0}$. Because all the stoichiometric compatibly classes are all translated from subspace $\mathscr{S}$, they share the same dimension. Any two functions $x_0(s)$ and $x_1(s)$ in $L_W$ we can find corresponding points $x_0(0)$ and $x_2$ in non-delayed system and the subtraction of the two points  is in the $Q_W\cap \mathscr{S}$, thus $\dim(Q_W\cap\mathcal{P}_{x_0})=s-1$. Therefore, $Q_W$ is a facet in stoichiometric compatibly class for mass action system without delay. Especially, when $W$ is a locking set, all the reactions are "locked" and all the points in the $Z_W$ are all equilibria. Then:
\begin{equation}
g(x_0(s))-g(x_1(s))=x_0(0)-x_1(0)
\end{equation}
Hence, if $x_0(s),~x_1(s)\in L_W$ are in the same stoichiometric compatibly class of delayed system, we have $x_0(0),~x_1(0)$ also in the same compatibly class of system without delay and vice versa.  In this case, it is obvious that $\dim(Q_W\cap\mathcal{P}_{x_0})=\dim(Z_W)=s-1$, namely, $Q_W\cap\mathcal{P}_{x_{0}}$ is a facet for system without delay.
$\Box$
\begin{lemma}\label{lm:facet}
For a DeCBMAS $(\mathcal{S,C,R},k,\tau)$ governed by \ref{eq:dde}, there does not exist an $\omega$-$limit$ set in the interior of $L_W$ if $L_W\cap \mathcal{D}_{\psi}$ is a facet of the stoichiometric compatibility class $\mathcal{D}_{\psi}$ in $(\mathcal{S,C,R},k,\tau)$, where $L_W$ is given by \ref{def:Lw} and $\psi$ is the initial data of the trajectory.
\end{lemma}
\noindent{\textbf{Proof.}}
We just need consider $L_W$ where $W$ is a semilocking set. Besides, without loss of generality, let $W=\{X_{d+1},...,X_{n}\}$. 
For delayed system, if there exists a trajectory $x^{\psi}(t)$ which trends to the interior of $L_W$, the $\omega$-$limit$ set can only be a connected equilibria set---$\omega(\psi)\subset L_W$ and for any $\epsilon>0$ we can find an $\epsilon$-neighbourhood of $\omega(\psi)$ and a large enough $t_0$ such that the trajectory will always stay in this $\epsilon$-neighbourhood when $t>t_0$. If this is not true, the trajectory will go out of the neighbourhood infinitely. Then an infinite function sequence $x^{\psi}_{t_n}(s),~s\in [-\tau_{max},0]~and~t_n\rightarrow \infty$ can be chosen with $dist(x^{\psi}_{t_n}, \omega(\phi))>\epsilon$ for all $t_n$. And all the trajectories in complex balanced system are bounded, so from Bolzano-Weierstrass theorem, there must exist an $\omega$-$limit$ point that out of the $\epsilon$-neighbourhood. But this is contradiction to the connectivity of $\omega$-$limit$ set. Thus when $t>t_0$, the trajectory $x(t)$ stays in $\epsilon$-neighbourhood for any $\epsilon>0$. As we know all the $\omega$-$limit$ functions are constant equilibria, in this case for any $t_k$ large enough, there exist $w_k\in \omega(\psi)$ and an $\epsilon_1$ such that:
$$dist(x^{\psi}_{t_k},~w_k)=\sup_{s\in [-\tau_{max},0]}\vert x^{\psi}_{t_k}(s)-w_k\vert<\epsilon_1/2.$$  
In this case, we can obtain that:
\begin{equation}\label{eq:tau}
 dist(x^{\psi}_{t_k}(0),~x^{\psi}_{t_k}(-\tau_{max}))< dist(x^{\psi}_{t_k}(0),w_k)+ dist(w_k,~x^{\psi}_{t_k}(-\tau_{max}))<\epsilon_1.    
\end{equation}
From \cite{Anderson2010} we know for zero-delayed system we can find a minimal complex $yl$(namely, $yl_{j}\leq y_{j}$ for all $j\in W$ and all $y\in L_{l}$ ) in each linkage class $L_l$. $yl$ is the dominate complex in all complexes in $L_l$. In other words, for any linkage class $L_l$ the relation between 0 and the following equation
\begin{equation*}
  \sum_{y_{.i}\rightarrow y'_{.i}\in L_l}k_i(x(t))^{y_{\cdot i}}(y'_{\cdot i}-y_{\cdot i})
\end{equation*}
is determined by $(y'_{\cdot i}-y_{\cdot i})$. For complex balanced system, each result complex is also reactant complex. So $yl'_{j}-yl_{j}\geq0$ for each $j\in W$ and all $l$. So we have
   \begin{equation}\label{eq:del}
   kl(x(t))^{yl}(yl'_{j}-yl_{j}) =kl(x(t))^{yl}yl'_{j}(1-\frac{yl_{j}}{yl'_{j}})>0.
\end{equation}
for some $j\in W$.
Then we consider the dynamics system with delay. 
\begin{align*}
    \dot{x}_j^{\psi}(t)&=\sum_{y_{\cdot i}\rightarrow y'_{\cdot i}\in \mathcal{R}}k_{i}\left(x(t-\tau_i)\right)^{y_{\cdot i}}y'_{ji}-\sum_{y_{\cdot i}\rightarrow y'_{\cdot i}\in \mathcal{R}}k_i\left(x(t)\right)^{y_{\cdot i}}y_{j i},~t\geq 0\\
    &=\sum_{l}\sum_{R_i\in L_l}\left(k_i(x(t-\tau_i))^{y_{\cdot i}}y'_{ji}-k_i(x(t))^{y_{\cdot i}}y_{ji}\right)\\
    &=\sum_{l}\sum_{\substack{R_i\in L_l\\y'_{ji}\neq 0}}k_i(x(t))^{y_{\cdot i}}y'_{ji}\left(\frac{x(t-\tau_i)^{y_{.i}}}{x(t)^{y_{.i}}}-\frac{y_{ji}}{y'_{ji}}\right)-\sum_{l}\sum_{\substack{R_i\in L_l\\y'_{ji}= 0}}k_i(x(t))^{y_{.i}}y_{ji}
 \end{align*}
The same as the non-delayed system, $yl$ is the minimal complex in each $L_l$. Thus above equation have the same sign with the following equation when $t$ is large enough:
\begin{equation}\label{eq:yl}
   \sum_{L_l}kl(x(t))^{yl}yl'_{j}(\frac{(x(t-\tau))^{yl}}{(x(t))^{yl}}-\frac{yl_{j}}{yl'_{j}})
\end{equation}
  where $kl$ is the reaction rate of $yl\rightarrow yl'$. For each $L_{l}$,
because $yl_{j}<yl'_{j}$, $\frac{yl_{j}}{yl'_{j}}<1$ is a constant. From \ref{eq:tau} we can choose a $t_1$ such that:
\begin{equation}\label{eq:ddel}
  \biggl{\vert} \frac{(x(t-\tau))^{yl}}{(x(t))^{yl}}-1\biggr{\vert}<\biggl{\vert} 1-\frac{yl_{j}}{yl'_{j}}\biggr{\vert}
\end{equation}
for all $t>t_1$.
We can derive $\dot{x}^{\psi}_{j}(t)>0$ from \ref{eq:yl}, \ref{eq:del} and \ref{eq:ddel} when $t$ large enough.  So the $\dot{x}^{\psi}_j(t)>0$ when $x^{\psi}(t)$ trends to its $\omega$-$limit$ set. But this is a contradiction to each function $\phi$ in $\omega$-$limit$ set with $\phi_j=0$ where $j\in W$. So our assumption is wrong. Any trajectory can not go to the interior of $L_W$.
$\Box$

We then get one of our main results in this paper.
\begin{theorem}\label{thm:facet}
Give a DeCBMAS $(\mathcal{S,C,R},k,\tau)$ of \ref{eq:dde}, for any semilocking set $W\subseteq\mathcal{S}$, if $L_W\cap \mathcal{D}_{\psi}$ defined by \ref{def:Lw} is either empty or a facet of the stoichiometric subspace $\mathscr{S}$, then this  DeCBMAS is persistent. 
\end{theorem}
\noindent{\textbf{Proof.}}
Supposing $W$ is a semilocking set and $W\subset W_1\subset\mathcal{S}$, we will proof $L_{W_1}$ can not exist $\omega$-limit points by using reduction to absurdity. So we suppose that $L_{W_1}$ is non-empty and $W_1$ is a semilocking set. Then $L_{W_1}$ must be a facet of stoichiometric subspace. And without loss of generality, let $W=\{X_1,~,...,X_{d}\}$ and $W_1=\{X_1,~...,X_{d+m}\}$. We assume $\{w_1,~...,w_{s-1},~v\}$ spans the stoichiometric subspace $\mathscr{S}$ and $\{w_1,...,w_{s-1}\}$ spans $\mathscr{S}\vert_{W_1}$. Then for any two $\psi_1, \psi_2\in \bar{\mathscr{C}}_{+}$, $$\left(g(\psi_1)-g(\psi_2)\right)\vert_{W^{c}}=qv\vert_{W_1}.$$ $v\vert_{W_1}$ is a basis spans $\mathscr{S}\vert_{W_1}$. and $q$ is a real number. $v\vert_{W_{1}}=(v_1,...,v_{d+m})$ and there does not exist species whose concentration never change. So all $v_1,...,v_{d+m}$ is not zero. But consider $\psi_1\in L_{W_1}$ and $\psi_2\in L_{W}$, we obtain:
$$(g(\psi_1)-g(\psi_2))\vert_{W_{1}}=\mathbf{g}=(\underbrace{0,...,~0}_{m},g_{m+1},...,g_{m+d})$$ with $g_{m+j},~j=1,...,d$ is nonzero. This obviously can not be expressed by $v\vert W_{1}$. So $L_{W_1}$ can not be a facet. Then $L_{W_1}$ does not exist $\omega$-$limit$ point. So we obtain if a trajectory trends to a facet $L_W$, there can not exist $\omega$-limit points in $L_{W_1}$ where $W\subset W_1$. Then combining with  $\ref{thm:wl}$ and $\ref{lm:facet}$ we can obtain the network is persistent when $Z_W$ is a facet or empty for all semilocking set $W$.
$\Box$
 
\subsection{DeCBMASs with $\dim{Z_W}=0$ for any semilocking set $W\subseteq \mathcal{S}$}
We consider another class of DeCBMASs that have persistence, identified by $\dim{Z_W}=0$. The result recurs that for complex balanced MASs presented by Anderson \cite{Anderson2008}.
 \begin{lemma}\label{lm:locking}
For a DeCBMAS $(\mathcal{S,C,R},k,\tau)$ modeled by \ref{eq:dde}, if $W\subseteq\mathcal{S}$ is a semilocking set and $L_W\cap \mathcal{D}_{\psi}$ is a vertex, i.e., $\dim{Z_W}=0$, then $W$ must be a locking set. Here, $L_W$ and $Z_W$ are defined accorrding to \ref{def:Lw} and \ref{eq:g}, respectively.
 \end{lemma}
 \noindent{\textbf{Proof.}}
  We will prove by contradiction. We assume that $W$ is a semilocking set but not a locking set. Without loss of generality, let $W=\{X_{n-n'+1},...,X_{n}\}$. Because $W$ is not a locking set, we can find the subnetwork $\mathcal{N_W}=(\mathcal{S_W,C_W,R_W})$ with 
 $$\mathcal{R}_W=\{y_{\cdot i}\vert_{W^{c}}\rightarrow y'_{\cdot i}\vert_{W^{c}}, \mathrm{supp}~y\cap W=\emptyset, \mathrm{supp}~y'\cap W=\emptyset\}.$$
 The set $\mathcal{S}_W$ and set $\mathcal{C}_W$ are the species and complexes appear in the $\mathcal{R}_W$. The definition of network tells us that reactions $y_{\cdot i}\rightarrow y'_{\cdot i}\in \mathcal{R}_W$ must satisfy $(y'_{\cdot i}-y_{\cdot i})\vert_{W^{c}}\neq 0$.
 So $\dim(\mathscr{S}\vert_{W^{c}})>0$. So we can choose two functions $\tilde{\psi}_1(s)$ and $\tilde{\psi}_2(s)$ in $\mathcal{N}_W$, such that:
 \begin{equation}\label{eq:0}
     g(\tilde{\psi}_1)-g(\tilde{\psi}_2)\neq \mathbbold{0}_{n-n'}\in \mathscr{S}_W.
 \end{equation}
 where $\mathscr{S}_W$ is the stoichiometric subspace of network $\mathcal{N}_W$.
 Then let $\psi_1=(\tilde{\psi}_1,\mathbbold{0}_{n'})^{\top}$ and $\psi_2=(\tilde{\psi}_2,\mathbbold{0}_{n'})^{\top}$, we can see $\psi_1$ and $\psi_2$ are all in $L_W\cap\mathcal{D}_{\psi}$. So we have the following equation:
 \begin{equation*} 
g(\psi_1)-g(\psi_2)=(g(\tilde{\psi}_1)-g(\tilde{\psi}_2), \mathbbold{0}_{n'})\in Z_W.
 \end{equation*}
 Above equation holds as that $g(\psi_1)-g(\psi_2)$ can be expressed by the linear combination of $$\{y'-y\vert~\mathrm{supp}~y\cap W=\emptyset,\mathrm{supp}~y'\cap W =\emptyset\}\subset \mathscr{S}.$$ 
 Above set doesn't equal to zero from \ref{eq:0}. So it is obviously contradict to $\dim Z_W=0$. $W$ must be a locking set. 
 $\Box$
 
 By using the chain method (See details in \ref{cm}), we can obtain an approximating system of the original delayed system. The following useful lemma researches the approximating system and will help us prove our main results.
 
 \begin{lemma}\label{lm:ad}
 For a complex balanced approximating system \ref{eq:cade}, if the point $(z^{\top},v^{\top})^{\top}=(\phi^{\top},\mathbbold{0}^{\top})^{\top}$ is an $\omega$-limit point and $W=\{Z_j\vert \phi_j=0\}$ is a locking set, we can find another point that has the same support set with $(z^{\top},v^{\top})^{\top}$.
 \end{lemma}
\noindent{\textbf{Proof.}}
 For arbitrary $\epsilon$, $T_{\epsilon}$ is the time that the trajectory $(z(t)^{\top},v(t)^{\top})^{\top}$ enter into $\epsilon$-neighbourhood of $(\phi^{\top},\mathbbold{0}^{\top})^{\top}$. $T_{\frac{\epsilon}{2}}$ is the time that trajectory enter into $\frac{\epsilon}{2}$-neighbourhood of $(\phi^{\top}, \mathbbold{0}^{\top})^{\top}$. The trajectory at time $T_{\epsilon}$ and $T_{\frac{\epsilon}{2}}$ are
 $\left(z(T_{\epsilon})^{\top},v^{\top}(T_{\epsilon})\right)^{\top},$ $\left(z(T_{\frac{\epsilon}{2}})^{\top},v^{\top}(T_{\frac{\epsilon}{2}})\right)^{\top}$, respectively.
 The Lyapunov function of approximating system is as follows: 
 \begin{equation*}
     H((z(t)^{\top},v(t)^{\top})^{\top})=\sum_{j=1}^{n}h_j(z_j)+\sum_{i=r-r'+1}^{r}\sum_{m=1}^{N}h_{im}(v_{im})
 \end{equation*}
 where $h_j(z_j)=\bar{z}_j-z_j-z_j\ln{\frac{\bar{z}_j}{z_j}}$, $h_{im}(v_{im})=\bar{v}_{im}-v_{im}-v_{im}\ln{\frac{\bar{v}_{im}}{v_{im}}}$ and $(\bar{z}^{\top},\bar{v}^{\top})^{\top}$ is the positive equilibirum in this stoichiometric compatibility class.
From the dissipation of function $H$, we have:
\begin{equation*}
    H\left((z(T_{\frac{\epsilon}{2}})^{\top},v(T_{\frac{\epsilon}{2}})^{\top})^{\top}\right)-H\left((z(T_{\epsilon})^{\top},v(T_{\epsilon})^{\top})^{\top}\right)\leq 0
\end{equation*}
Assuming $W=\{Z_1, ...,Z_d\}$ with $d<n$, above equation can be rewritten as:
\begin{align*}
      &H\left((z(T_{\frac{\epsilon}{2}})^{\top},v(T_{\frac{\epsilon}{2}})^{\top})^{\top}\right)-H\left((z(T_{\epsilon})^{\top},v(T_{\epsilon})^{\top})^{\top}\right)\\
       =&\sum^{d}_{j=1}\left(h_{j}\left(z_{j}(T_{\frac{\epsilon}{2}})\right)-h_{j}(z_{j}(T_{\epsilon}))\right)+\sum^{n}_{k=d+1}\left(h_{k}(z_{k}(T_{\frac{\epsilon}{2}}))-h_{k}(z_{k}(T_{\epsilon}))\right)\\
      &+\sum_{i=r-r'+1}^{r}\sum^{N}_{m=1}(h_{im}(v_{im}(T_{\frac{\epsilon}{2}}))-h_{im}(v_{im}(T_{\epsilon})))\\
      =&\underbrace{\sum^{d}_{j=1}\left(\ln(\tilde{z}_{j})-\ln(\bar{z}_{j})\right)(z_{j}(T_{\frac{\epsilon}{2}})-z_j(T_{\epsilon}))}_{a}+\sum^{n}_{k=d+1}\left(\ln(\tilde{z}_{k}))-\ln(\bar{z}_{k})\right)(z_{k}(T_{\frac{\epsilon}{2}})-z_{k}(T_{\epsilon})))\\
     &+\underbrace{\sum_{\mathrm{supp}~y_{\cdot i}\cap W\neq \emptyset}\sum^{N}_{m=1}(h_{im}(v_{im}(T_{\frac{\epsilon}{2}}))-h_{im}(v_{im}(T_{\epsilon})))}_{b}\leq 0.\\
\end{align*}
The last equation holds from differential mean value theorem and $\tilde{z}_{j}\in (z_j(T_{\frac{\epsilon}{2}}),$
$z_j(T_{\epsilon}))$. 
We can choose $\epsilon$ small enough such that for all the $Z_j\in W$, there exist $\tilde{z}_j<\epsilon<\bar{z}_j$. And in this case, $a\geq 0$. Denote $z_k(T_{\frac{\epsilon}{2}})-z_k(T_{\epsilon})$ as $\Delta z_k$ and we have:
\begin{equation*}\small
    \left\vert a \right\vert
    \leq \left\vert\sum^{n}_{k=d+1}\left(\ln(\tilde{z}_{k}))-\ln(\bar{z}_{k})\right)\Delta z_k
    +\sum_{\mathrm{supp}~y_{\cdot i}\cap W\neq \emptyset}\sum^{N}_{m=1}(h_{im}(v_{im}(T_{\frac{\epsilon}{2}}))-h_{im}(v_{im}(T_{\epsilon})))\right\vert\\
\end{equation*}
For each $j\leq d$, if let $\vert\ln(\bar{z}_{j})\vert=M_{j}$, $\left\vert\ln(\tilde{z}_{j})-\ln(\bar{z}_{j})\right\vert\geq \vert\ln{\epsilon}\vert-M_{j}$. When $\epsilon\rightarrow 0$, $\vert\ln(\epsilon)\vert-M_{j}$ trends to infinity. And $h_{im}(v_{im})$ trends to infinity only when $v_{im}$ goes to infinity. But from the boundedness of complex balanced and $W$ is a locking set, we can obtain:
\begin{equation}\label{eq:bb}
    \vert\Delta z_{j}\vert\leq \frac{1}{\vert\ln(\epsilon)\vert-M_j}\left(\sum_{k=d+1}^{n}c_{k}\vert\Delta z_{k}\vert+M\right), \forall j\leq d
\end{equation}
where $c_k=\ln(\tilde{z}_{k})-\ln(\bar{z}_{k})$ and $M$ is the upper bound of part $b$. Further we assume  $\Delta_{max}=\sup_{k=d+1,\cdots,n}\{\vert\Delta z_{k}\vert\}$ , $\delta(\epsilon)=\sup_{j\in \{1,\cdots,d\}}(1/\vert\ln(\epsilon)-M_j\vert)$ and $C=\sum_{k=d+1}^{n}c_{k}$. Then for $j=1,\cdots,d$, \ref{eq:bb} can be written as:
\begin{equation*}
    \vert\Delta z_j\vert\leq \delta(\epsilon)(C\Delta_{max}+M).
\end{equation*}
Now it is the time to consider the vector:
$$(z^1(\epsilon)^{\top},v^1(\epsilon)^{\top})\triangleq\frac{\Delta (z^{\top},v^{\top})}{\Delta_{max}}=\left(\frac{z(T_{\frac{\epsilon}{2}})}{\Delta_{max}}^{\top},\frac{v(T_{\frac{\epsilon}{2}})}{\Delta_{max}}^{\top}\right)-\left(\frac{z(T_{\epsilon})}{\Delta_{max}}^{\top},\frac{v(T_{\epsilon})}{\Delta_{max}}^{\top}\right)\in \mathscr{S}'$$
where $\mathscr{S}'$ is the stoichiometric subspace of approximating system.
For $j=1,\cdots,d$,
\begin{equation}\label{eq:id}
    \vert z^1_{j}(\epsilon)\vert\leq\delta(\epsilon)D
\end{equation}
 where $D$ is a positive constant. Besides, as there must exist some $k\in\{d+1,\cdots,n\}$ such that $z_k=\Delta_{max}$, then $1\leq\vert (z^1(\epsilon),v^1(\epsilon))\vert\leq M_1$ for some positive constant $M_1$. Because $\epsilon$ is small enough and arbitrary, we can choose a decreasing sequence $\{\epsilon_n\}$ with $\epsilon_n\rightarrow 0$. Then from above analysis, we can also obtain a sequence of vectors $\{(z^1(\epsilon_{n}),v^{1}(\epsilon_{n})\}$. And each vector has one $\delta(\epsilon_n)$ such that $\delta(\epsilon_n)\rightarrow 0$ when $\epsilon_{n}\rightarrow 0$. From \ref{eq:id} we can further derive that $z^{1}_{j}(\epsilon_{n}),~j\in \{1,\cdots,d\}$ and $v^{1}_{im}$ where $i\in\{r-r'+1,\cdots,r\}, m\in\{1,\cdots,N\}$ also trends to zero as $\epsilon_{n}\rightarrow 0$. $v^1_{im}\rightarrow 0$ follows that $W$ is a locking set and each $v^1_{im}$ can be expressed as:$C_0z^{1}_{j}$ for some constant $C_0$ and some $j\in\{1,\cdots,d\}$.
For each $n$, vector $(z^1(\epsilon_{n}),v^{1}(\epsilon_{n}))^{\top}$ is contained in compact space 
$$O=\{(z^1,v^1)^{\top}:1\leq \vert (z^1,v^{1})^{\top}\vert\leq M_1\}\cap \mathscr{S}'.$$
So there exist a convergent sub-sequence $\{(z^1(\epsilon_{n_{j}})^{\top},v^{1}(\epsilon_{n_{j}})^{\top}\}$ and it converges to the point $((z^0)^{\top},(v^0)^{\top})^{\top}\in\mathscr{C}\cap\mathscr{S}'$ when $j\rightarrow +\infty$. From above analysis, $\vert z^0\vert>1$ and have the following form:
\begin{equation*}
  ((z^0)^{\top},(v^0)^{\top})^{\top}=(\mathbbold{0}_{d},
  z^{0}_{d+1},\cdots,z^{0}_{n},\mathbbold{0}_{r'N})^{\top}\in\mathscr{S}' 
\end{equation*}
Because the stoichiometric subspace of the approximating system $\mathscr{S}'=\mathscr{S}\oplus\mathscr{S}_{1,i}\oplus\mathscr{S}_{2,i}$, where $\mathscr{S}_{1,i}$ have no impact on $z^0$ and because $v_{i,n+1}=0$ for all $i\in\{r-r'+1,\cdots,r\}$, $\mathscr{S}_{2,i}$ also have no impact on $z^0$. Thus $z^0$ is in $\mathscr{S}$. 
    And then for any $\alpha>0$, $((z^2)^{\top},(v^2)^{\top})^{\top}=(\phi^{\top},\mathbbold{0}_{r'N}^{\top})^{\top}+\alpha((z^0)^{\top},(v^0)^{\top})^{\top}$ in the same stoichiometric compatibility class with $(\phi^{\top},\mathbbold{0}_{r'N}^{\top})^{\top}$ also have the same support set with $(z^{\top},v^{\top})^{\top}=(\phi^{\top},\mathbbold{0}_{r'N}^{\top})^{\top}$.
 $\Box$
 
 We thus obtain another main result of this paper.
 
\begin{theorem}\label{thm:sv}
 For a DeCBMAS $(\mathcal{S,C,R},k,\tau)$ given by \ref{eq:dde}, if $L_W\cap \mathcal{D}_{\psi}$ is a vertex for each semilocking set $W\subseteq\mathcal{S}$, this system is persistent, where $L_W$ and $\mathcal{D}_{\psi}$ share the same meanings with those in \ref{lm:locking}.
 \end{theorem}
 \noindent{\textbf{Proof.}}
 Lemma \ref{lm:locking} reveals that if $L_W$ only have one element in each stoichiometric compatibly class $\mathcal{D}_{\psi}$, then $W$ must be a locking set. Without loss of generality, we assume that reactions $\{y_{\cdot i}\rightarrow y'_{\cdot i},~i=r-r'+1,~...,r\}$ in system $\mathcal{N}$ have time delay. Instead of considering the original system $\mathcal{N}=(\mathcal{S,~C,~R},~k,~\tau)$, we study the approximating system $\mathcal{N}'=(\mathcal{S',~C',~R'},~k')$ derived from adding a chain of first-order reactions to each delayed reactions in original system.
   And the dynamical equation of $\mathcal{N}'$ is \ref{eq:cade}. And when $\phi(s)$ is the unique vertex in $L_W\cap \mathcal{D}_{\psi}$, the corresponding point $(z^{\top},v^{\top})^{\top}$ in approximating system is:
 \begin{equation*}
     z^{\top}=\tilde{\phi}(s),
   v_{im}^{\top}=k_i\frac{\tau_i}{N}\tilde{z}^{y_{\cdot i}}=k_i\frac{\tau_i}{N}\tilde{\phi}^{y_i}=0.
 \end{equation*}
 Above equation is true from the approximation of \ref{app} and $\tilde{\phi}^{y_{\cdot i}}\equiv 0$ for all the points in $L_W$ when $W$ is a locking set. $\tilde{\phi}$ is an $\omega$-limit point, $(\tilde{\phi}^{\top},\mathbbold{0}^{\top})^{\top}$ is also an $\omega$-limit point in the approximating system. So from \ref{lm:ad}, we can find $(z^2,v^2)=(\tilde{\phi},(\mathbbold{0}_{r'N})^{\top})^{\top}$. In this way, when $N\rightarrow +\infty$, $(z^2,v^2)$ corresponding to a constant function $\tilde{\phi}(s)$ for all $s\in [-\tau_{max},0]$ in $L_W$. And we can compute $$g(\phi)-g(\tilde{\phi})=z^0\in\mathscr{S}.$$
    In other words, $\tilde{\phi}$ is also in $L_W\cap\mathcal{D}_{\theta}$. This is contradicting to the condition. So the unique point of intersection can not be $\omega$-$limit$ point. Thus system $\mathcal{N}$ is  persistent.
 $\Box$
 
Gathering up \ref{thm:facet},\ref{thm:sv}, we have a comprehensive result about the persistence of DeCBMASs. 
 \begin{theorem}\label{th:evf}
 Consider a DeCBMAS $(\mathcal{S,C,R},k,\tau)$ described by \ref{eq:dde}. Then the DeCBMAS is persistent if $L_W\cap\mathcal{D}_{\psi}$ is empty, a facet or a vertex of the stoichiometric compatibility class $\mathcal{D}_{\psi}$ where $L_W$ is given by \ref{def:Lw} and $\psi\in\mathscr{C}_{+}$ is an initial function for each semilocking set $W\subseteq\mathcal{S}$.
  \end{theorem}
  \noindent{\textbf{Proof.}}
 Through the proof of \ref{lm:facet}, there can not exists an $\omega$-limit point in the interior of $L_W$ if $L_W\cap\mathcal{D}_{\psi}$ is a facet. \ref{thm:sv} reveals that the unique function $\phi$ in $L_W\cap\mathcal{\psi}$ can not be an $\omega$-limit point. So this result is obvious.
  $\Box$
 \begin{corollary}
  All 2-dimensional DeCBMASs $(\mathcal{S,C,R},k,\tau)$ are persistent. 
 \end{corollary}
\noindent{\textbf{Proof.}}
For any 2-dimensional DeCBMAS $(\mathcal{S,C,R},k,\tau)$, the dimension of $Z_W$ defined by \ref{eq:g} is constrained by $\dim{Z_W}=0$ or $\dim{Z_W=1}$. Therefore, the result is true. 
$\Box$
 
 It has been proven that each stoichiometric compatibility class in a complex balanced MAS contains a unique positive equilibrium. Moreover, this equilibrium is locally asymptotically stable, and further globally asymptotically stable if persistence property is added \cite{Horn1972}. We exhibit this result to be also true for DeCBMASs.
 
 \begin{theorem}\label{th:gas}
 Consider a DeCBMAS $(\mathcal{S,C,R},k,\tau)$ of \ref{eq:dde}, and for each semi-locking set $W\subseteq\mathcal{S}$ with $L_W$ defined by \ref{def:Lw}, $L_W\cap\mathcal{D}_{\psi}$ is empty, a facet or a vertex of the  compatibility class $\mathcal{D}_{\psi}$, where $\psi\in\mathscr{C}_{+}$ is an initial function and $\bar{x}\in\mathbb{R}^n_{>0}$ is the corresponding positive equilibrium. Then $\bar{x}$ is globally asymptotically stable with respect to all initial functions in $\mathcal{D}^+_{\psi}$ given in \ref{eq:Dpsi}.
 \end{theorem}
  \noindent{\textbf{Proof.}}
 We can easily conclude this result from \ref{thm:wl} and \ref{th:evf}.
   $\Box$

\section{Some examples}\label{sec5}
In this section we will illustrate the previous two main results through two examples.

\newtheorem{example}{Example}
\begin{example}\label{ex:1}
Consider a delayed network $(\mathcal{S,C,R},k,\tau)$ with the reaction route
\begin{align*}
\xymatrix{2X_1\ar ^{\tau_1,~k_1~~~} [r] & 3X_{1}+X_{2} \ar ^{~\tau_2,~k_2} [d]\\
 & X_{1}+2X_2\ar ^{~\tau_3,~k_3~~} [lu]}
\end{align*}
It is easy to check that this network is of weakly reversible structure with the deficiency to be zero, and thus should be complex balanced regardless of the rate constant vector $k$ \cite{Horn1972}. Further, its stoichiometric subspace $\mathscr{S}$ is found to be $\mathbb{R}^2$, i.e.,
$$\mathscr{S}=\mathrm{span}\{(1,1)^{\top}, (-1,2)^{\top}\}.$$ 
By setting $k=(1,1,2)^{\top}$, the dynamics follows
\begin{equation}
\begin{split}
     \dot{x}=&\left(x_1(t-\tau_1)\right)^2\left(
     \begin{matrix}
     3\\
     1
     \end{matrix}
     \right)+\left(x_1(t-\tau_2)\right)^3x_2(t-\tau_2)\left(
     \begin{matrix}
     1\\
     2
     \end{matrix}
     \right)+x_1(t-\tau_3)\left(x_2(t-\tau_3)\right)^2\left(
     \begin{matrix}
     2\\
     0
     \end{matrix}
     \right)\\
     &-\left(x_1(t)\right)^2\left(
     \begin{matrix}
     2\\
     0
     \end{matrix}
     \right)-\left(x_1(t)\right)^3x_2(t)\left(
     \begin{matrix}
     3\\
     1
     \end{matrix}
     \right)-x_1(t)\left(x_2(t)\right)^2\left(
     \begin{matrix}
     1\\
     2
     \end{matrix}
     \right)
\end{split}
\end{equation}
This DeCBMAS has a unique positive equilibrium $\bar{x}=(1.49,0.95)^{\top}$.

The species set $\mathcal{S}$ of this network include two subsets to be the semilocking sets, respectively being $W_1=\{X_1\}$ and $W_2=\{X_1,X_2\}$. According to \ref{def:Lw}, we could  write
$$L_{W_1}=\{\psi\in\bar{\mathscr{C}}_{+}\vert \psi (s)=(0,\psi_2(s))^\top, s\in[-\tau_{\max},0],\psi_2(s)\neq 0\}~\text{and}~L_{W_2}=\{(0,0)^\top\}.$$
Further based on \ref{eq:g}, we could know that the representative element $v\in Z_{W_1}$ has the first entry to be 0 while the second entry to be undetermined. This means that any element $v\in Z_{W_1}$ can be linearly expressed by  the vector $(0,-1)^{\top}$, i.e., $\dim{Z_{W_1}}=1=\dim\mathscr{S}-1$, so $L_{W_1}\cap\mathcal{D}_\psi$ is a facet of $\mathcal{D}_\psi$. Also since $L_{W_2}\cap\mathcal{D}_\psi=(0,0)^\top$, i.e., only containing one element, $L_{W_2}\cap\mathcal{D}_\psi$ is a vertex of $\mathcal{D}_\psi$. According to \ref{th:evf}, this DeCBMAS is persistent.   

\begin{figure}
        \centering
        \begin{minipage}[c]{0.46\textwidth}
        \centering
          \includegraphics[width=1\textwidth]{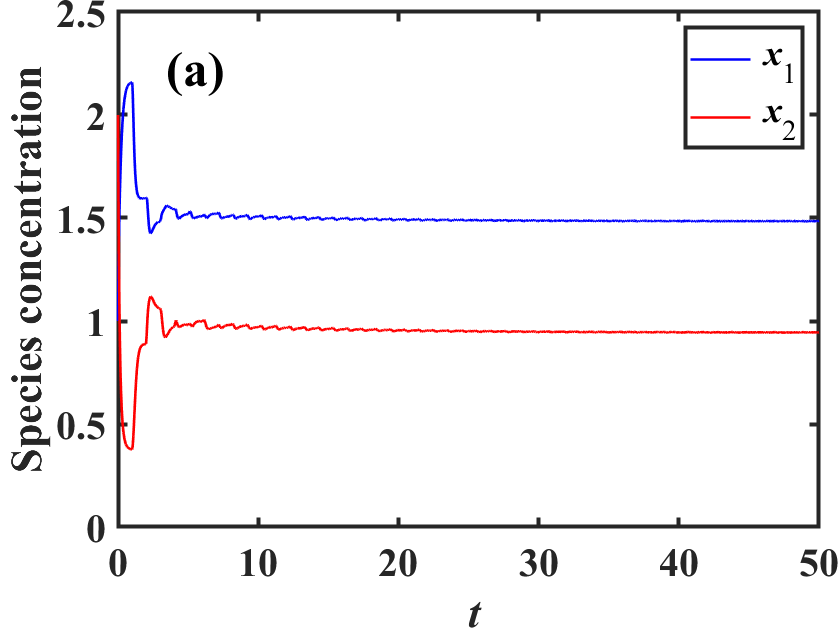}
          \end{minipage}
          \begin{minipage}[c]{0.46\textwidth}
          \centering
          \includegraphics[width=1\textwidth]{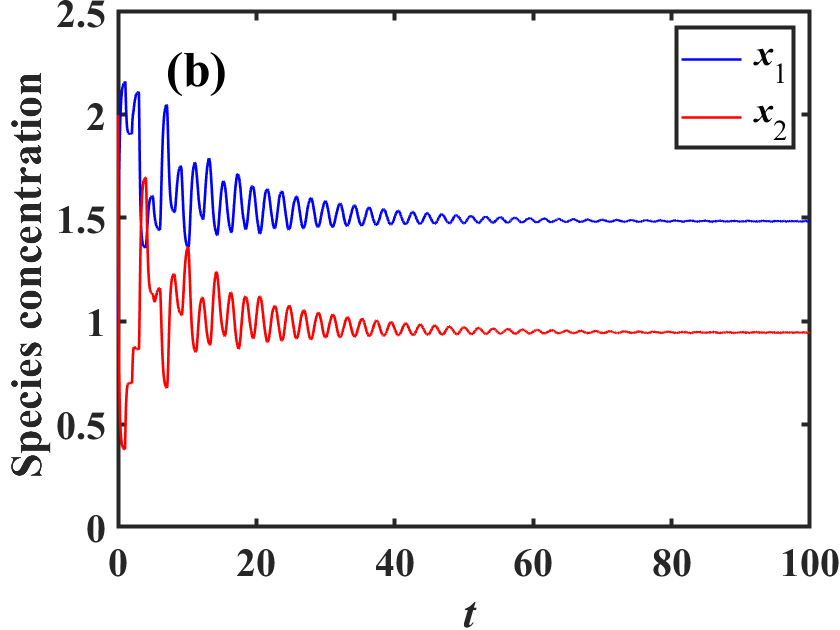}
        \end{minipage}
         \begin{minipage}[c]{0.46\textwidth}
          \centering
          \includegraphics[width=1\textwidth]{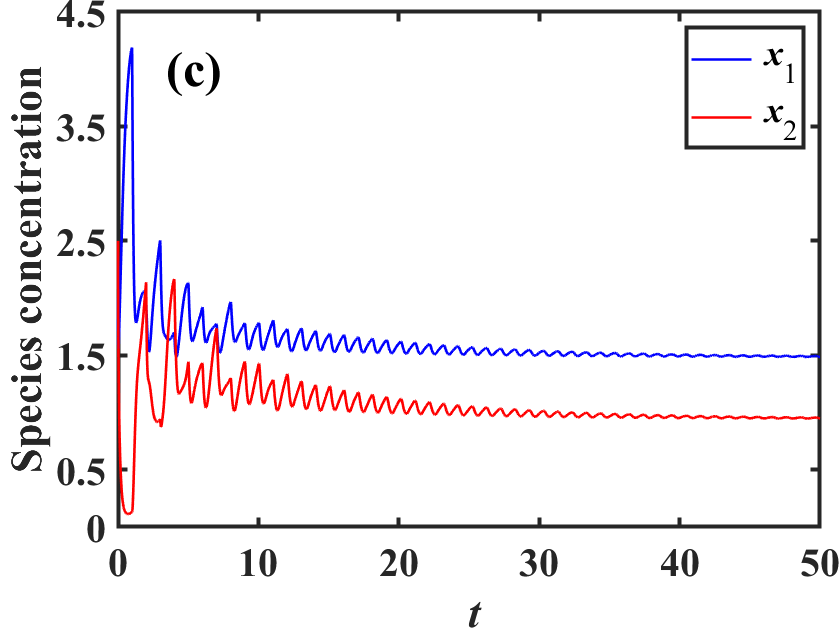}
        \end{minipage}
         \begin{minipage}[c]{0.46\textwidth}
          \centering
          \includegraphics[width=1\textwidth]{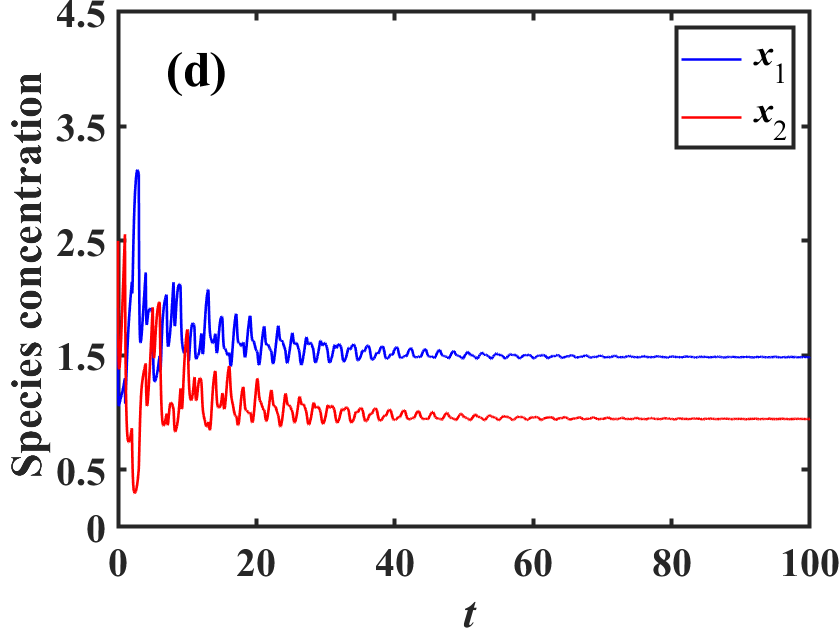}
        \end{minipage}
        \vspace{2em}
         \caption{The evolution behavior of \ref{ex:1} with different initial data and different time delay: (a) $\tau=(1,2,1),~\psi(s)=(1,2)$; (b) $\tau=(2,1,3),~\psi(s)=(1,2)$;  (c) $\tau=(1,2,1),~\psi(s)=(\sin(s)+2,\cos(s)+1.5)$; (d) $\tau=(2,1,3),~\psi(s)=(\sin(s)+2,\cos(s)+1.5)$.}\label{fig:1}
    \end{figure}

We exhibit the evolution behavior of this system at different time delays and different initial functions. Shown in \ref{fig:1} are the results. As can be seen, the trajectories will converge to $\bar{x}$ although they start from different initial functions and the system are assigned different time delay. 
\end{example}

\begin{example}\label{ex:2}
\begin{figure}
        \centering
        \begin{minipage}[c]{0.46\textwidth}
        \centering
          \includegraphics[width=1\textwidth]{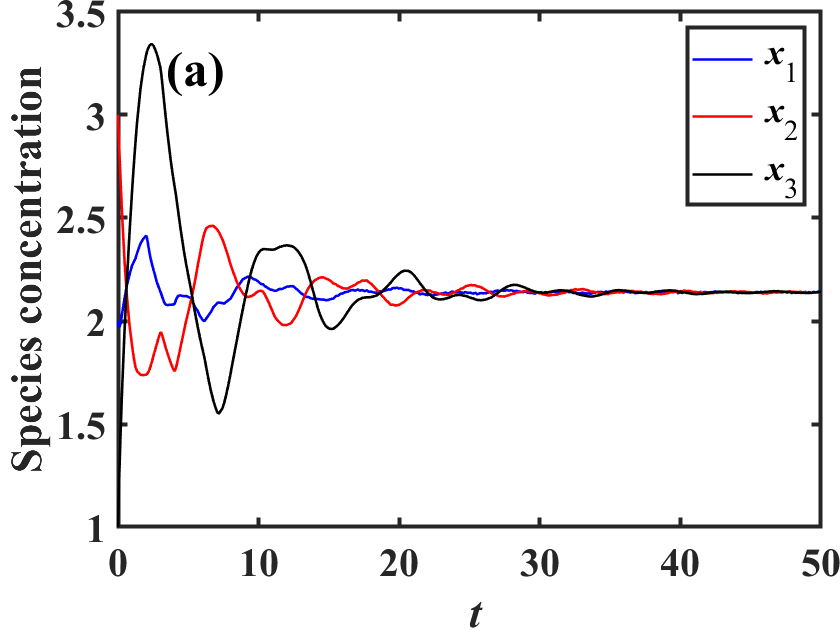}
          \end{minipage}
          \begin{minipage}[c]{0.46\textwidth}
          \centering
          \includegraphics[width=1\textwidth]{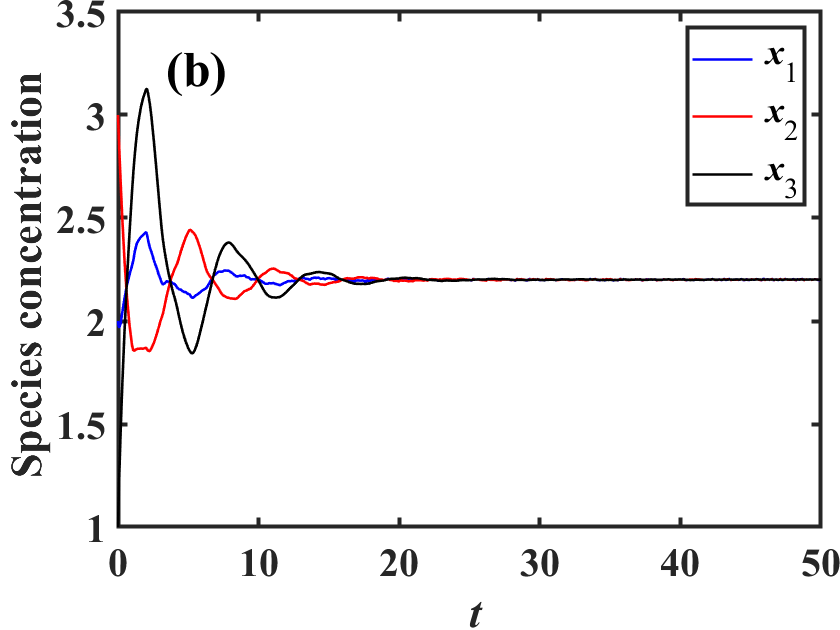}
        \end{minipage}
         \begin{minipage}[c]{0.46\textwidth}
          \centering
          \includegraphics[width=1\textwidth]{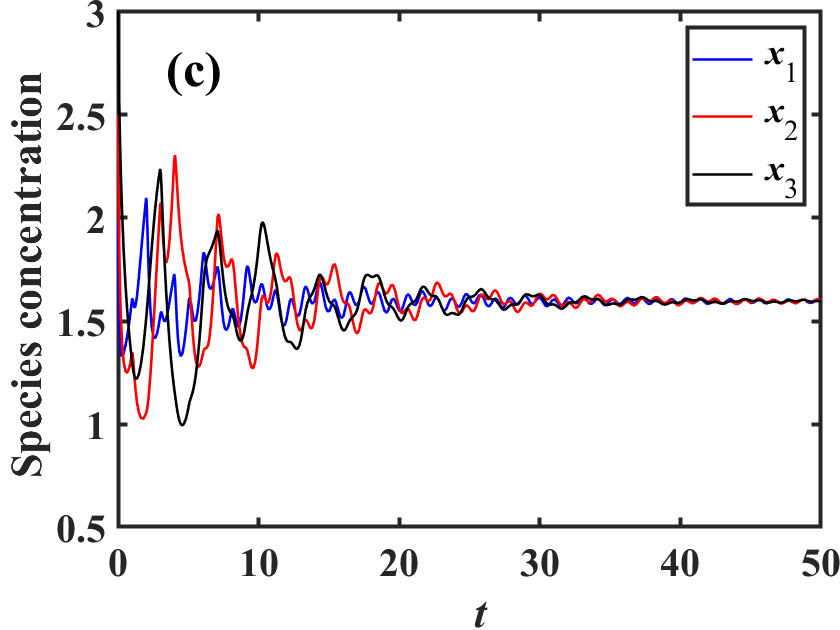}
        \end{minipage}
         \begin{minipage}[c]{0.46\textwidth}
          \centering
          \includegraphics[width=1\textwidth]{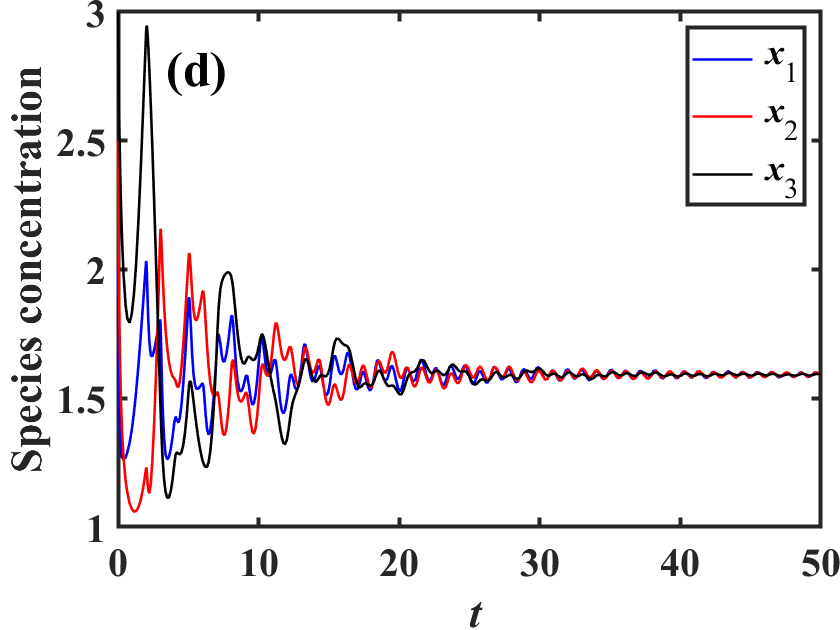}
        \end{minipage}
        \vspace{2em}
         \caption{The evolution behavior of \ref{ex:2} with different initial data and different time delay: (a) $\tau=(1,2,3,4),~\psi(s)=(2,3,1)$; (b) $\tau=(3,2,2,1),~\psi(s)=(2,3,1)$; (c) $\tau=(1,2,3,4),~\psi(s)=(\sin(s)+2,\cos(s)+1.5,2\sin(s)+3)$; (d) $\tau=(3,2,2,1),\psi(s)=(\sin(s)+2,\cos(s)+1.5,2\sin(s)+3)$.}\label{fig:2}
    \end{figure}

The second example is a DeCBMAS $(\mathcal{S,C,R},k,\tau)$:
\begin{equation*}
  \xymatrix{2X_1 \ar @{ -^{>}}^{\tau_1,k_1~}  @< 1pt> [r] & X_1+X_2  \ar @{ -^{>}}^{\tau_2,k_2~} @< 1pt> [l] \ar @{ -^{>}}^{~\tau_3,k_3~}  @< 1pt> [r] & X_{2}+X_3 \ar @{ -^{>}}^{~\tau_4,k_4~} @< 1pt> [l]},
 \end{equation*}
The stoichiometric subspace is $\mathscr{S}=\mathrm{span}\{(-1, 1, 0)^{\top},(-1,0,1)^{\top}\}$ and $\dim\mathscr{S}=2$. Also assume $k=(1,1,1,1)^{\top}$, then the dynamical equation takes
\begin{equation*}
\begin{split}
    \dot{x}(t)=&x_1^2(t-\tau_1)\left(
     \begin{matrix}
     1\\
     1\\
     0
     \end{matrix}
     \right)+x_1(t-\tau_2)x_2(t-\tau_2)\left(
     \begin{matrix}
     2\\
     0\\
     0
     \end{matrix}
     \right)+x_1(t-\tau_3)x_2(t-\tau_3)\left(
     \begin{matrix}
     0\\
     1\\
     1
     \end{matrix}
     \right)\\+&x_2(t-\tau_4)x_3(t-\tau_4)\left(
     \begin{matrix}
     1\\
     1\\
     0
     \end{matrix}
     \right)
     -x_1^2(t)\left(
     \begin{matrix}
     2\\
     0\\
     0
     \end{matrix}
     \right)-2x_1(t)x_2(t)\left(
     \begin{matrix}
     1\\
     1\\
     0
     \end{matrix}
     \right)-x_2(t)x_3(t)\left(
     \begin{matrix}
     0\\
     1\\
     1
     \end{matrix}
     \right).
\end{split}
\end{equation*}
The system can reach complex balancing no matter what $\tau$ is, and each complex balanced equilibrium $\bar{x}$ satisfies $\bar{x}_1=\bar{x}_2=\bar{x}_3$. 

There are three semilocking sets $W_1=\{X_1,X_2\}$, $W_2=\{X_1,X_3\}$ and $W_3=\{X_1,X_2,X_3\}$ for this network. We thus get 
\begin{align*}
 L_{W_1}=&\{\psi\in\bar{\mathscr{C}}_{+}\vert \psi (s)=(0,0,\psi_3(s))^\top, s\in[-\tau_{\max},0],\psi_3(s)>0\};\\
    L_{W_2}=&\{\psi\in\bar{\mathscr{C}}_{+}\vert \psi (s)=(0,\psi_2(s),0)^\top, s\in[-\tau_{\max},0],\psi_2(s)>0\};\\
   L_{W_3}=&\{(0,0,0)^\top\}.
\end{align*}
Clearly, $L_{W_1}\cap\mathcal{D}_{\psi}$, $L_{W_2}\cap\mathcal{D}_{\psi}$ and $L_{W_3}\cap\mathcal{D}_{\psi}$ are all vertexes of $\mathcal{D}_{\psi}$. From \ref{th:evf}, it suggests that this DeCBMAS is persistent. \ref{fig:2} also exhibits the persistence of the system starting from different initial function and different time delay.  
\end{example}

\section{Conclusions}\label{sec6}
This paper has tackled an issue of analyzing the persistence of DeCBMASs. For a DeCBMAS, if $L_W\cap\mathcal{D}_{\psi}$ is empty, a facet or a vertex of $\mathcal{D}_{\psi}$ with respect to each semilocking set $W$ in this network, then the DeCBMAS is persistent. These results recur those for diagnosing the persistence of complex balanced MASs presented by Anderson et al. \cite{Anderson2008,Anderson2010}. Namely, a complex balanced MAS has persistence if $Q_W\cap\mathcal{P}_{x_0}$ is empty, a facet of $\mathcal{P}_{x_0}$ or discrete. Further, we prove the DeCBMASs admitting the proposed conditions for persistence are also globally asymptotically stable at each positive equilibrium relative to the corresponding positive stoichiometric compatibility class. This result also recurs that for complex balanced MASs, i.e., they are globally asymptotically stable at each positive equilibrium relative to its positive stoichiometric compatibility class if they are persistent. 

We thus suppose whether there are time delays in reactions that will not affect the persistence of reaction network systems. The following research may focus on other MASs that have been proved to have persistence, like strong endotactic network \cite{G2014}, to check persistence when time delays in reactions are modeled.      

\begin{appendix}
\section{Chain method}\label{cm} In this appendix, we recall the chain method for the chemical reaction kinetics \cite{Repin1965}.
 
Consider a delayed system with the following equation
 \begin{equation*}
     \dot{x}(t)=f(x(t))+f_1(x(t-\tau))J,~\forall~ t\geq0,
 \end{equation*}
 where $x(t)\in \mathbb{R}^{n}$ is the state vector,  $f:\mathbb{R}^{n}\rightarrow \mathbb{R}^{n}$ and $f_1:\mathbb{R}^{n}\rightarrow \mathbb{R}$ are continuous functions, $\tau>0$, and $J$ is a constant vector. Let $x(t)=\psi(t)$ for $-\tau\leq t \leq 0$ be a continuous initial function. Then we use a set of ODEs with a ``$N$" new state variables chain to approximate the original delayed system. The ODEs have the following form
 \begin{equation}\label{eq:ade}
     \begin{split}
         \dot{z}(t)&=f(z(t))+\frac{N}{\tau}v_{N}(t)J,\\
     \dot{v_1}(t)&=f_1(z(t))-\frac{N}{\tau}v_1(t),\\
     \dot{v_j}(t)&=\frac{N}{\tau}v_{j-1}(t)-\frac{N}{\tau}v_{j}(t),~2\leq j\leq N,
     \end{split}
     \end{equation}
where $z(t)\approx x(t)\in\mathbb{R}^{n}$, $v_j~(1\leq j \leq N)$ are the added chain, $z(0)=\psi(0)$ and $v_j(0)=\int^{-(j-1)\frac{\tau}{N}}_{-j\frac{\tau}{N}}f_1(\psi(s))ds$. Repin \cite{Repin1965} revealed that if the initial function of the delayed system is sufficiently smooth, the solution of above ODEs converges uniformly to the solution of the original delayed model on any finite time interval $[0,T]$ when $N$ goes to infinity. For a DeMAS $(\mathcal{S,C,R},k,\tau)$, we assume that the last $r'$ reactions have delays, i.e., $\tau_i=0$ for $i=1,...,r-r'$ and $\tau_i>0$ for the rest of $i$. In this case, the delayed dynamical equation can be written as
 \begin{equation*}
\dot{x}(t)=\sum^{r-r'}_{i=1}k_i(x(t))^{y_{\cdot i}}\left[y'_{\cdot i}-y_{\cdot i}\right]+\sum^{r}_{i=r-r'+1}k_i\left[(x(t-\tau_i))^{y_{\cdot i}}y'_{\cdot i}-(x(t))^{y_{\cdot i}}y_{\cdot i}\right].
 \end{equation*}
By applying \ref{eq:ade}, we write out the approximating system to be
 \begin{equation}
 \begin{split}
\label{eq:cade}
     \dot{z}(t)=\sum^{r-r'}_{i=1}k_i(z(t))^{y_{\cdot i}}[y'_{\cdot i}-y_{\cdot i}&]+\sum^{r}_{i=r-r'+1}\frac{N}{\tau_i}v_{iN}(t)y'_{\cdot i}-\sum_{i=r-r'+1}^{r}k_i(x(t))^{y_{\cdot i}}y_{\cdot i},\\
     \dot{v}_{i1}(t)&=k_i(z(t))^{y_{\cdot i}}-\frac{N}{\tau_{\cdot i}}v_{i1}(t),\\
     \dot{v}_{ij}(t)=\frac{N}{\tau_i}&v_{i,j-1}(t)-\frac{N}{\tau_i}v_{j}(t),~2\leq i\leq N,
 \end{split}
 \end{equation}
 with the initial points
 \begin{equation*}
     z(0)=\psi(0),~\text{and}~v_{ij}(0)=k_j\int^{-(j-1)\frac{\tau_i}{N}}_{-j\frac{\tau_i}{N}}(\psi(s))^{y_{\cdot i}}ds.
 \end{equation*}
 In the above dynamical equation, we take each time delay $\tau_k$ as a chain of first-order intermediate reactions with each reaction rate to be $\frac{N}{\tau_i}$, and further let $\hat{y}_{\cdot i}=[y^{\top}_{\cdot i}, \mathbbold{0}^{\top}]^{\top}$ and $\hat{y}'_{\cdot i}=[y'^{\top}_{\cdot i},\mathbbold{0}^{\top}]^{\top}$ for $i=1,...,r$. Then the stoichiometric subspace, denoted by $\mathscr{S'}$, can be expressed by :
 $$\mathscr{S'}=\mathscr{S}\oplus \mathscr{S}'_{1,i}\oplus\mathscr{S}'_{2,i},~i=r-r'+1,...,r,$$
 where
 \begin{equation*}
     \begin{split}
     \mathscr{S}&=\mathrm{span}\left(\{\hat{y}'_{\cdot i}-\hat{y}_{\cdot i}\vert i=1,...,r\}\right),\\
     \mathscr{S}'_{1,i}&=\mathrm{span}\left(\{\textbf{e}_{i,j+1}-\textbf{e}_{i,j}\vert j=(n+1),...,(n+N-1)\}\right),\\
     \mathscr{S}'_{2,j}&=\mathrm{span}\left(\{\textbf{e}_{i,n+1}-\hat{y}_{i}\}\right).
 \end{split}
 \end{equation*}
Every point $\hat{x}(t)$ in the approximating system can be written as
 $$\hat{x}(t)=(z(t)^{\top},v(t)^{\top})^{\top}=(z(t)^{\top},v_{r-r'+1,1}(t),...,v_{r-r'+1,N}(t),...,~v_{r,1}(t),...,~v_{r,N}(t))^{\top}.$$
Liptk et al. \cite{G2018} tell us that by using the properties of the chain method as $N\rightarrow \infty$, we have the following approximations
 \begin{equation}\label{app}
 \begin{split}
       z(t)&\approx x(t),\\
        v_{ij}(t)\approx k_i\int^{t-(j-1)\frac{\tau_{i}}{N}}_{t-j\frac{\tau_{i}}{N}}&(x(s))^{y_{\cdot i}}ds\approx k_i\frac{\tau_i}{N}\left(x(t-j\frac{\tau_i}{N})\right)^{y_{\cdot i}}.
 \end{split}
 \end{equation}
\end{appendix}

\end{document}